%% file: ex_article_arvix_nosupp.tex
\documentclass[hidelinks,onefignum,onetabnum]{siamart251216}

\input{ex_shared}
\usepackage{xcolor} 
\newcommand{\pendingbox}[1]{\fbox{\parbox{0.9\textwidth}{\centering\footnotesize\textcolor{red}{Pending: \detokenize{#1}}}}}
\newcommand{\maybegraphics}[2][]{\IfFileExists{#2}{\includegraphics[#1]{#2}}{\pendingbox{#2}}}
\newcommand{\maybetable}[1]{\IfFileExists{#1}{\resizebox{\textwidth}{!}{\input{#1}}}{\pendingbox{#1}}}

\ifpdf
\hypersetup{
  pdftitle={Adaptive Blockwise Tensor Trains for Stochastic Finite-Volume Methods with Uncertain Geometry},
  pdfauthor={WP Gammel, ST Walton, SA Andrews, S Tokareva, WT Taitano}
}
\fi

\begin{document}

\maketitle

%
%
\begin{abstract}
Stochastic finite-volume (SFV) methods extend conservative, shock-capturing discretizations to hyperbolic systems with uncertain inputs by partitioning the parameter domain into finite-volume cells. However, their memory and computational cost grow exponentially with the number of uncertain parameters. Tensor-train (TT) formats can remove this exponential dependence, provided the ranks of the discrete solution remain small throughout the simulation. This requirement is frequently violated when the uncertainty is geometric. Mapping each realization onto a common reference domain introduces parameter-dependent metric factors, so that shocks and contacts propagate at different reference speeds and come to occupy different cells across realizations. The resulting rank growth is  spatially localized---a structure that neither a single global TT nor a separate TT for each physical cell is able to exploit. In this work, we propose an adaptive blockwise TT representation in which the SFV solution is partitioned along its physical mode into contiguous blocks, each stored as an independent TT for every conserved variable. Blocks are split near discontinuities and in regions of large variance or rank, and merged where the solution is smooth, with separate refinement and coarsening thresholds preventing oscillation of the partition. Unlike block-structured adaptive mesh refinement, the mesh itself is never changed; only the grouping of cells into trains adapts. We embed this representation in a TT-based SFV solver combining WENO3 reconstruction in the physical and stochastic directions, positivity-preserving limiters, and strong-stability-preserving Runge--Kutta time integration, with nonlinear steps evaluated via cross approximation. A mean-conservative rounding that adds no rank, together with single-valued interface fluxes, conserves the parameter-mean mass and energy to round-off for any partition. The approach is demonstrated on a Sod shock-tube problem with uncertain wall positions. 
\end{abstract}



%
%
\section{Introduction}

Uncertainty quantification (UQ) carries the variability of model inputs, such as material properties, initial and boundary data, and the shape of the computational domain, through PDE-based simulations to the outputs of interest, supplying moments, distributions, and confidence intervals that make simulation results comparable to experimental data and usable in risk-informed decisions \cite{delia2020quiet,smith2013uq}. Hyperbolic conservation laws, such as those that describe compressible gas dynamics, shallow-water flows, and multi-material hydrodynamics, are a particularly challenging setting for UQ. In this setting, perturbations of the inputs shift the arrival time and strength of shocks \cite{abgrall2017uq,poette2009uq}, and
discontinuities can form in finite time even from smooth data. This creates a solution that is non-smooth in both the stochastic variables and physical space \cite{schwab2013shock,abgrall2017uq}. Scalable UQ methods for these systems must estimate statistics reliably despite low parametric regularity \cite{mishra2012sparse}, remain stable across shocks and contacts \cite{leveque2002fvm,kusch2020filtered}, and have a cost that grows slowly with the number of uncertain parameters \cite{abgrall2017sfv,bachmayr2016tensor,mishra2012mlmc}.

Existing methods are usually distinguished by how strongly they modify the underlying deterministic solver \cite{geraci2016weakly}. Non-intrusive methods treat the solver as a black box (e.g., Monte Carlo sampling and its multilevel extensions \cite{giles2008mlmc,mishra2012mlmc,mishra2012sparse,hoel2012amlmc,giles2015distribution,eigel2016amlmc,elfverson2016failure,krumscheid2018mlmc}). While they are robust, they converge algebraically in the number of samples. Intrusive methods reformulate the governing equations in the stochastic variables (e.g., by Galerkin projection onto polynomial chaos \cite{knio2006pc} or multiresolution bases \cite{herty2024multiresolution}). However, for discontinuous solutions, global polynomial bases converge slowly and produce Gibbs-type oscillations, and the projected system can lose hyperbolicity \cite{pettersson2015pc,kusch2020filtered}. The ``weakly intrusive'' stochastic finite-volume (SFV) method \cite{tokareva2014sfv,abgrall2017sfv,tokareva2024gas,harmon2025adaptive,harmon2024active} sits between these two classes. It discretizes the parameter domain with finite-volume cells and evolves averages over products of physical and stochastic cells, so that an existing finite-volume code can be reused with modest changes. Discrete conservation and hyperbolicity carry over by construction, and high-order accuracy in the parameters follows from reconstructions local to each stochastic cell, without a global smoothness assumption \cite{abgrall2017sfv,harmon2025adaptive}. However, with $M$ cells per parameter and $d$ parameters, each stochastic cell stores $M^d$ values of every conserved variable, so the cost of SFV grows exponentially with the stochastic dimension \cite{abgrall2014te,abgrall2017sfv,walton2025ttsfv}.

Low-rank tensor formats reduce this cost whenever the discrete solution is approximately separable across its modes, as shown for ordinary and partial differential equations \cite{dolgov2012parabolic,kazeev2012laplace,khoromskij2011qtt,dolgov2015cme,kazeev2017qttfem,kazeev2018qttfem,dolgov2019conservation,manzini2023maxwell,truong2023neutron} and for UQ \cite{dolgov2015pce,zhang2015anova}; see the surveys \cite{grasedyck2013survey,kolda2009tensor}. In the tensor-train (TT) format \cite{oseledets2011tt,oseledets2010cross}, memory and basic arithmetic grow only linearly with the number of modes for small ranks, linear operators act on the cores without approximation, SVD-based rounding lowers ranks with controlled error, and nonlinear pointwise functions can be sampled by cross approximation without assembling the full array. Writing the SFV unknowns in TT form thus compresses the stochastic dimensions while keeping the conservative, shock-capturing update intact, provided the ranks stay small throughout the simulation.

Geometric uncertainty, such as uncertainty in the domain length, layer thicknesses, and interface positions, makes this condition hard to satisfy. Prior work \cite{gammel2026ttsfv} maps each realization onto a common reference domain, so that a single stationary mesh serves all realizations and the geometry appears only through Jacobian factors in the flux. Since the Jacobian depends on the parameters, so do the reference wave speeds, and shocks and contacts come to occupy different cells for different realizations. If the SFV solution is arranged as a matrix with one row per physical cell and one column per stochastic cell, a discontinuity whose position varies over $n_{\mathrm{w}}$ physical cells produces $n_{\mathrm{w}}$ shifted step profiles, none expressible in terms of the others, so the numerical rank at any tolerance below the jump height is of order $n_{\mathrm{w}}$. This rank growth is confined to where the parameter-dependent waves pass, but neither classical TT layout exploits that. If the entire domain is stored in one TT \cite{walton2025ttsfv}, a single leading rank links every physical cell to the stochastic modes; since TT arithmetic scales cubically in the ranks and Hadamard products multiply them, the few cells containing discontinuities dictate the cost everywhere. If every physical cell carries its own TT \cite{dubois2026ttsfv, gammel2026ttsfv}, rank growth remains local, but the parametric structure common to neighboring cells is lost and the parametric cores are replicated across the mesh.

In this work, we introduce an adaptive blockwise TT representation for SFV that balances these two extremes, in the spirit of hierarchical and block low-rank matrix formats \cite{hackbusch1999hmatrix,amestoy2015blr}. The SFV solution tensor is sliced along the physical mode into blocks of contiguous cells, each stored, for each conserved variable, as a separate TT whose leading mode is the block-local cell index. Neighboring cells within a block share parametric cores, a discontinuity in one block does not affect the ranks of the others, and the leading rank of a block of $n_b$ cells cannot exceed $\min\{n_b, P^d\}$, where $P^d$ is the number of stochastic cells attached to each physical cell. Blocks are halved where indicators based on jumps in the mean state, the density variance, or the TT rank exceed refinement thresholds, and adjacent blocks are merged where all indicators fall below smaller coarsening thresholds, so that the partition does not alternate between steps. A split restricts the physical core and a merge concatenates two trains, each followed by rounding, and a merge is kept only if rounding perturbs the density TT by less than a set tolerance. The mesh and the discretization stay fixed: adaptation alters how the solution is compressed, not how finely it is resolved, unlike adaptive mesh refinement \cite{berger1984amr,berger1989amr} and moving-mesh methods \cite{hirt1974ale,lohner1996ale}.

Block-local TTs have recently been combined with block-structured adaptive mesh refinement (AMR) for thermal radiation transport \cite{gorodetsky2026blocktt}, where the space--angle intensity in each AMR block is stored as its own TT. Both works rest on the observation that decomposing the domain into blocks confines non-separable structure to the blocks where it occurs, but they differ in four respects. First, the blocks in \cite{gorodetsky2026blocktt} are AMR blocks of equal size and refinement changes the resolution, whereas here the mesh is fixed, the blocks vary in size, and only the grouping of cells into trains adapts. Second, the transport operator and block transfers in \cite{gorodetsky2026blocktt} are linear, whereas our scheme, with WENO weights, positivity limiters, and approximate Riemann solvers, is nonlinear throughout, and every nonlinear step is evaluated by cross approximation. Third, conservation in \cite{gorodetsky2026blocktt} holds to within the TT rounding tolerance, whereas our mean-conservative rounding, which adds no rank, together with single-valued interface fluxes conserves the parameter-mean mass and energy to round-off for any partition. Fourth, the high-rank structure arises there from angular transport, whereas here it arises from the parameter-dependent metric factors of geometric uncertainty.

We test the approach on a Sod shock tube whose two reflecting walls are independent uniform random variables, while the diaphragm is held fixed in physical space. The map to the unit interval aligns the walls but places the diaphragm at a parameter-dependent reference location, so the initial jump spans part of the reference domain and the parameter-dependent wave speeds spread the shock and contact further, while the flow elsewhere is close to deterministic. The main contributions of this work are the following.
\begin{itemize}
\item An adaptive blockwise TT representation for SFV on a fixed mesh, whose partition into blocks of variable size is adapted by exact split and merge operations on the physical cores, driven by shock, uncertainty, and rank indicators with separate refinement and coarsening thresholds and a test on the rounding error of each merge.
\item A fully nonlinear TT-based SFV solver for the Euler equations on the partition, with WENO3 reconstruction in space and parameters \cite{jiang1996weno,shu1998eno}, positivity-preserving limiters \cite{zhang2010positivity,zhang2011positivity}, Rusanov or HLL fluxes \cite{rusanov1962,harten1983hll}, and SSP-RK3 time stepping \cite{shu1988ssp,gottlieb2001ssp}. Blocks are coupled by a two-cell halo exchange, and every nonlinear step is evaluated by DMRG-type cross approximation with maxvol pivoting \cite{oseledets2010cross,savostyanov2011dmrg,goreinov2010maxvol}.
\item  A mean-conservative TT rounding that restores the parameter-mean totals without increasing the ranks. Together with single-valued interface fluxes, it conserves the parameter-mean mass and energy to round-off for any partition, independently of the TT tolerance.
\item A comparison of the blockwise, global, and per-cell TT layouts on an uncertain-geometry Sod problem, covering the accuracy of the statistics, both in the reference coordinate and mapped back to physical space, conservation, storage, and runtime.
\end{itemize}

Section~\ref{sec:FV_scheme} describes the SFV discretization of the reference-domain Euler equations, Section~\ref{sec:tt} its realization in TT arithmetic, including the mean-conservative rounding, and Section~\ref{sec:blockwise} the adaptive blockwise representation. Section~\ref{sec:results} reports results for the uncertain-geometry Sod problem, Section~\ref{sec:discussion} discusses them, and Section~\ref{sec:conclusions} concludes.

%
%
\section{Finite Volume Scheme under Geometric Uncertainty}
\label{sec:FV_scheme}
Following \cite{gammel2026ttsfv}, we build geometric uncertainty, i.e., uncertainty in the location and length of the physical domain, into the discretization from the outset. We restrict the discussion to one physical dimension.

\subsection{Governing equations on an uncertain interval}
\label{sec:geo_unc}
Let $\mathbf{y} = (y_1, y_2) \in \Gamma \subset \mathbb{R}^{2}$ be two geometric parameters that set the endpoints of the physical interval,
\begin{equation}
D_{x}(\mathbf{y}) = \big[x_L(\mathbf{y}),\, x_R(\mathbf{y})\big], \quad x_L(\mathbf{y}) = y_1, \quad x_R(\mathbf{y}) = y_2,
\label{eqn: geo_un_int}
\end{equation}
of length $s(\mathbf{y}) := x_R(\mathbf{y}) - x_L(\mathbf{y})$. Each $\mathbf{y} \in \Gamma$ describes one geometry, a \emph{realization}. For every realization we seek a solution $\mathbf{U}(t,x;\mathbf{y})$ of the compressible Euler equations
\begin{equation}
    \partial_{t} \mathbf{U} + \partial_{x} \mathbf{F}(\mathbf{U}) = \mathbf{0}; \ (t,x) \in (0, T] \times D_{x}(\mathbf{y}), \qquad
    \mathbf{U} = \begin{pmatrix} \rho \\ \rho u \\ \rho e \end{pmatrix}, \quad
    \mathbf{F}(\mathbf{U}) = \begin{pmatrix} \rho u \\ \rho u^{2} + p \\ u(\rho e  + p ) \end{pmatrix},
\label{eqn:comp_euler_1D}
\end{equation}
with density $\rho$, velocity $u$, total energy per unit volume $\rho e$, and the ideal-gas pressure $p = (\gamma - 1)(\rho e - \frac{1}{2}\rho u^{2})$, with adiabatic index $\gamma = 1.4$.

\subsection{Mapping to a fixed reference interval}
\label{sec:ref_map}
Discretizing Eqn.~\ref{eqn:comp_euler_1D} directly would require a different mesh for every realization, so that statistics in a given cell would require interpolation between incompatible meshes. We therefore map every realization to the reference interval $D_{\xi} = [0,1]$ by the affine map $x(\xi; \mathbf{y}) := x_L(\mathbf{y}) + s(\mathbf{y})\,\xi$, whose Jacobian $\partial \xi/\partial x = 1/s(\mathbf{y}) =: \alpha(\mathbf{y})$ is constant in space. By the chain rule, $\partial_x \mathbf{F} = \alpha(\mathbf{y})\,\partial_{\xi} \mathbf{F}$, and Eqn.~\ref{eqn:comp_euler_1D} becomes
\begin{equation}
\partial_t \mathbf{U}(t, \xi; \mathbf{y}) + \alpha(\mathbf{y})\,\partial_{\xi} \mathbf{F}(\mathbf{U}(t, \xi; \mathbf{y})) = \mathbf{0}, \quad (t, \xi) \in (0, T] \times (0,1).
\label{eqn:cons_law_ref}
\end{equation}
Since $\alpha$ does not vary with $\xi$, Eqn.~\ref{eqn:cons_law_ref} remains in conservative form. The factor $\alpha$ rescales every wave speed in reference coordinates: on a short domain ($s < 1$) a wave appears faster in $\xi$, and on a long domain the reverse holds.

\subsection{Mapping statistics back to physical space}
\label{sec:stat_map}
Statistics computed at a fixed $\xi$ compare the realizations at the same relative spot in its domain, i.e., at different physical points $x(\xi; \mathbf{y})$. Conversely, a feature that is fixed in physical space, such as the diaphragm of our test problem (see: Section~\ref{sec:ic_bc}), will occupy a range of different reference positions (Eqn.~\ref{eqn:xi_disc}). As a result, statistics computed in $\xi$ mix two kinds of uncertainty: the 'true' uncertainty of the flow, and uncertainty introduced by the mapping itself. To separate the two, we report statistics computed at both fixed $\xi$ and $x$. Statistics computed at fixed $x$ are obtained from the inverse map
\begin{equation}
\xi(x;\mathbf{y}) = \frac{x - x_L(\mathbf{y})}{s(\mathbf{y})}.
\label{eqn:inv_map}
\end{equation}
Since wall positions vary between realizations, a physical point $x$ is not necessarily inside every domain. The realizations whose domains contain $x$ form the set
\begin{equation}
\Gamma_x := \big\{\mathbf{y} \in \Gamma :\; x_L(\mathbf{y}) \le x \le x_R(\mathbf{y})\big\}.
\label{eqn:cover_set}
\end{equation}
Thus, statistics at $x$ are conditional on $\mathbf{y}\in \Gamma_{x}$. Let $\mu$ denote the probability measure of $\mathbf{y}$, which in this work is uniform on $\Gamma$. The mean of a quantity $g \in \{\rho, u, p\}$ at $(t, x)$
is then
\begin{equation}
\mathbb{E}_x[g](t,x) := \frac{1}{\mu(\Gamma_x)} \int_{\Gamma_x}
g\big(t, \xi(x;\mathbf{y}); \mathbf{y}\big)\, d\mu(\mathbf{y}),
\label{eqn:phys_mean}
\end{equation}
that is, the average of $g$ over only those realizations that contain $x$. The
variance and percentiles are defined in the same way, from the conditional
distribution of $g(t, \xi(x;\mathbf{y}); \mathbf{y})$ on $\Gamma_x$. We call
$\mu(\Gamma_x)$, the \emph{coverage} of $x$. For the walls of Section~\ref{sec:sfv_grids}, the
coverage equals one for $\delta_L \le x \le 1 - \delta_R$, where every
realization contains $x$. It decreases to zero towards $x = -\delta_L$ and
$x = 1 + \delta_R$. Because the map is
affine and invertible for each realization, Eqns.~\ref{eqn:inv_map}--\ref{eqn:phys_mean} only require post-processing of the reference-domain solution, and the solver itself is unchanged. The discrete
evaluation is described in Section~\ref{sec:results_setup}.

\subsection{Discretization Scheme}
\label{sec:sfv_grids}
We partition the reference interval into $N$ uniform cells $K_{i}^{\xi} = [\xi_{i - \frac{1}{2}}, \xi_{i + \frac{1}{2}}]$ with centers $\xi_i$ and spacing $\Delta \xi = 1/N$. In the SFV method, the parameters are treated as additional coordinates and discretized by finite-volume cells as well \cite{tokareva2014sfv,abgrall2017sfv}. Labeling the two parameter directions by the walls $L$ and $R$, the parameter domain is $\Gamma = \Gamma_L \times \Gamma_R$ with $\Gamma_L = [-\delta_L, \delta_L]$ and $\Gamma_R = [1 - \delta_R, 1 + \delta_R]$, where $\delta_L$ and $\delta_R$ are the maximal wall displacements. Each range is divided into $M$ uniform cells of width $h_L = 2\delta_L/M$ and $h_R = 2\delta_R/M$, and two \emph{padding} cells are appended on each side, giving $M + 4$ cells per direction with centers $y_{L,j_L}$ and $y_{R,j_R}$. A parameter cell is the rectangle $Y_{\mathbf{j}} = Y_{L,j_L} \times Y_{R,j_R}$ with area $|Y_{\mathbf{j}}| = h_L h_R$ and multi-index $\mathbf{j} = (j_L, j_R)$. Padding cells are evolved like interior cells; they only give interior cells near the edge of $\Gamma$ a complete reconstruction stencil and are excluded from all reported statistics. Since $s(\mathbf{y})$ can become small at the Gauss nodes of padding cells when $\delta_L$ and $\delta_R$ are large, we evaluate the geometric factor everywhere as
\begin{equation}
\alpha(\mathbf{y}) = \frac{1}{\max\{s(\mathbf{y}),\, s_{\min}\}}, \qquad s_{\min} = 0.4.
\label{eqn:alpha_clamp}
\end{equation}
For geometries with $s(\mathbf{y}) \ge s_{\min}$, Eqn.~\ref{eqn:alpha_clamp} reducies to the standard Jacobian.

The unknowns are the averages of the conserved state over space--parameter cells,
\begin{equation}
\bar{\mathbf{U}}_{i,\mathbf{j}}(t) := \frac{1}{\Delta\xi\,|Y_{\mathbf{j}}|} \int_{Y_{\mathbf{j}}} \int_{K_i^{\xi}} \mathbf{U}(t, \xi; \mathbf{y})\, d\xi \, d\mathbf{y}.
\label{eqn:ref_avg_cons}
\end{equation}
Integrating Eqn.~\ref{eqn:cons_law_ref} over $K_i^{\xi} \times Y_{\mathbf{j}}$, the averages change only through the $\alpha$-weighted fluxes across the two faces in $\xi$, averaged over $Y_{\mathbf{j}}$. There is no flux between parameter cells, which interact only through the reconstruction. Since $\alpha$ varies across $Y_{\mathbf{j}}$, it must be kept inside the parameter integral. Replacing that integral by Gauss quadrature with nodes $\mathbf{y}_{\mathbf{j},\mathbf{q}}$ and weights $w_{\mathbf{q}}$ (Section~\ref{sec:tt_sfv}), and the exact interface flux by a numerical flux $\hat{\mathbf{F}}(\cdot,\cdot)$ \cite{leveque2002fvm,vanleer1979muscl}, gives the semi-discrete SFV scheme
\begin{align}
\mathbf{G}_{i+\frac{1}{2},\mathbf{j}}(t) &:= \sum_{\mathbf{q}} w_{\mathbf{q}}\, \alpha(\mathbf{y}_{\mathbf{j},\mathbf{q}})\, \hat{\mathbf{F}}\Big(\mathbf{U}^{L}_{i+\frac{1}{2}}(t;\mathbf{y}_{\mathbf{j},\mathbf{q}}),\, \mathbf{U}^{R}_{i+\frac{1}{2}}(t;\mathbf{y}_{\mathbf{j},\mathbf{q}})\Big), \label{eqn:sem_dis_geounc}\\
\frac{d}{dt}\bar{\mathbf{U}}_{i,\mathbf{j}}(t) &= \big(\mathcal{L}(\bar{\mathbf{U}})\big)_{i,\mathbf{j}} := -\frac{1}{\Delta\xi}\Big(\mathbf{G}_{i+\frac{1}{2},\mathbf{j}}(t) - \mathbf{G}_{i-\frac{1}{2},\mathbf{j}}(t)\Big), \label{eqn:fv_fluxdiff}
\end{align}
where $\mathbf{U}^{L}_{i+\frac12}$ and $\mathbf{U}^{R}_{i+\frac12}$ are the states immediately to the left and right of $\xi_{i+\frac12}$. For a fixed $\mathbf{y}$, this is the familiar Godunov-type update on a physical mesh of width $s(\mathbf{y})\Delta\xi$. The interface states are obtained in two WENO3 stages, first in $\xi$ (Section~\ref{sec:weno_phys}), giving face values that are still parameter-cell averages, and then in $\mathbf{y}$ (Section~\ref{sec:weno_stoch}), giving values at the Gauss nodes. Both act on each conserved variable separately and are followed by a positivity-preserving limiter.

\subsection{WENO3 reconstruction in the reference coordinate}
\label{sec:weno_phys}
We use the standard third-order WENO reconstruction \cite{jiang1996weno,shu1998eno,shu1999weno}, which blends two linear candidates with weights that favor the smoother stencil and thereby avoids spurious oscillations near shocks \cite{sweby1984flux,shu1987tvb}. Fix $\mathbf{j}$ and let $\bar{w}_i$ denote one component of $\bar{\mathbf{U}}_{i,\mathbf{j}}$. With the smoothness indicators $\beta_{0,i} = (\bar{w}_i - \bar{w}_{i-1})^2$ and $\beta_{1,i} = (\bar{w}_{i+1} - \bar{w}_i)^2$ of the stencils $\{i-1,i\}$ and $\{i,i+1\}$, the values at the right and left faces of cell $i$ are
\begin{equation}
\begin{split}
w^-_{i+\frac12} &= \chi^-_{0,i}\Big(\tfrac{3}{2}\bar{w}_i - \tfrac{1}{2}\bar{w}_{i-1}\Big) + \chi^-_{1,i}\Big(\tfrac{1}{2}\bar{w}_i + \tfrac{1}{2}\bar{w}_{i+1}\Big), \\
w^+_{i-\frac12} &= \chi^+_{0,i}\Big(\tfrac{1}{2}\bar{w}_{i-1} + \tfrac{1}{2}\bar{w}_i\Big) + \chi^+_{1,i}\Big(\tfrac{3}{2}\bar{w}_i - \tfrac{1}{2}\bar{w}_{i+1}\Big),
\end{split}
\label{eqn:right_recon}
\end{equation}
with the nonlinear weights
\begin{equation}
\chi^{\pm}_{k,i} = \frac{\tilde{\chi}^{\pm}_{k,i}}{\tilde{\chi}^{\pm}_{0,i}+\tilde{\chi}^{\pm}_{1,i}}, \qquad \tilde{\chi}^{\pm}_{k,i} = \frac{\mu^{\pm}_k}{(\beta_{k,i} + \varepsilon)^2} \qquad k = 0,1.
\label{eqn:weno_weights}
\end{equation}
Here, $\mu^{\pm}_k$ are linear weights, which are chosen so that, if used unmodified, the blend reproduces the most accurate
value obtainable from all three cells, $\{i-1, i, i+1\}$.

Evaluating the polynomials at the interfaces in Eqn.~\eqref{eqn:right_recon} shows the reconstructions do not depend on the grid spacing. Thus, the reconstruction is unaffected by the mapping, and geometry enters the spatial discretization only through the Jacobian. 

The two states at interface $\xi_{i+\frac{1}{2}}$ are taken from the two cells that share it:
\begin{equation}
\begin{split}
\bar{\mathbf{U}}^{L}_{i+\frac{1}{2},\mathbf{j}} &= \bar{\mathbf{U}}^{-}_{i+\frac{1}{2},\mathbf{j}} \;\;(\textrm{right face of cell } i), \\ 
\bar{\mathbf{U}}^{R}_{i+\frac{1}{2},\mathbf{j}} &= \bar{\mathbf{U}}^{+}_{i+\frac{1}{2},\mathbf{j}} \;\;(\textrm{left face of cell } i+1),
\end{split}
\label{eqn:recon_LR_cons}
\end{equation}
where the overbar is retained because these traces are still averages over the parameter cell $Y_{\mathbf{j}}$.

WENO3 does not guarantee positive face values near strong jumps. We therefore apply the scaling limiter of \cite{zhang2010positivity,zhang2011positivity} to the density and total energy. Letting $w^{\min}_i$ $= \min\{w^{+}_{i-\frac{1}{2}}, w^{-}_{i+\frac{1}{2}}\}$ be the smaller face value of cell $i$ and $\varepsilon_w = 10^{-8}$ a small floor, we define
\begin{equation}
\theta_i = \begin{cases}
\min\Big\{1, \dfrac{\bar{w}_i - \varepsilon_w}{\bar{w}_i - w^{\min}_i}\Big\}, & \bar{w}_i > w^{\min}_i, \\[6pt]
1, & \textrm{otherwise},
\end{cases}
\qquad
w^{\pm} \leftarrow \bar{w}_i + \theta_i\,\big(w^{\pm} - \bar{w}_i\big),
\label{eqn:pos_lim_phys}
\end{equation}
with $\theta_i$ clipped to $[0,1]$. The momentum, which may take either sign, is not limited.

\subsection{WENO3 reconstruction in the parameter coordinates}
\label{sec:weno_stoch}
In the test problem defined in  Section~\ref{sec:ic_bc}, the initial discontinuity sits at the fixed physical location $x = 1/2$, so its reference position
\begin{equation}
\xi_d(\mathbf{y}) = \frac{1/2 - x_L(\mathbf{y})}{s(\mathbf{y})}
\label{eqn:xi_disc}
\end{equation}
depends on the geometry, as do the positions of the waves it generates. In a reference cell near a wave, the solution therefore jumps as a function of $\mathbf{y}$. A linear reconstruction across such a jump would generate Gibbs-type oscillations at the Gauss nodes, contaminating both the fluxes and the computed statistics. We therefore apply WENO3 in the parameter directions as well. Consider one direction $m$ with all other indices fixed, and let $\bar{v}_j$ denote one component of a face state from Section~\ref{sec:weno_phys}, averaged over parameter cell $j$. The reconstruction is evaluated at the local coordinate $\zeta = (y_m - y_{m,j})/h_m \in [-\tfrac{1}{2}, \tfrac{1}{2}]$ of the two Gauss points $\zeta_{1,2} = \pm 1/(2\sqrt{3})$,
\begin{equation}
v_j(\zeta) = \omega_0(\zeta)\big[\bar{v}_j + \zeta(\bar{v}_j - \bar{v}_{j-1})\big] + \omega_1(\zeta)\big[\bar{v}_j + \zeta(\bar{v}_{j+1} - \bar{v}_j)\big],
\label{eqn:weno_stoch}
\end{equation}
where the nonlinear weights $\omega_k(\zeta)$ are formed as in Eqn.~\ref{eqn:weno_weights} from the smoothness indicators $\beta_0 = (\bar{v}_j - \bar{v}_{j-1})^2$, $\beta_1 = (\bar{v}_{j+1} - \bar{v}_j)^2$, and the position-dependent ideal weights
\begin{equation}
\mu_0(\zeta) = \frac{\zeta - \zeta^2 + \frac{1}{12}}{2\zeta}, \qquad \mu_1(\zeta) = \frac{\zeta^2 + \zeta - \frac{1}{12}}{2\zeta}.
\label{eqn:weno_stoch_lin}
\end{equation}
For smooth data, the blend reproduces the quadratic whose averages over cells $j-1$, $j$, $j+1$ match the data. At the Gauss nodes, $\zeta_q^2 = 1/12$, such that $\mu_0 = \mu_1 = \frac{1}{2}$. There, the curvature term of the quadratic vanishes, and the smooth-data reconstruction reduces to the central-slope value $\bar{v}_j + \zeta_q(\bar{v}_{j+1} - \bar{v}_{j-1})/2$. Both ideal weights are therefore positive at every point where the reconstruction is evaluated. 

At the outermost padding cells, the reconstruction falls back to the constant value $\bar{v}_j$. Since $\mathbf{y}$ has two components, the reconstruction is applied one direction at a time. The two directions are swept in turn, the second producing the four point values $\mathbf{U}^{L/R}_{i+\frac{1}{2}}(\mathbf{y}_{\mathbf{j},\mathbf{q}})$ required by Eqn.~\ref{eqn:sem_dis_geounc}. 

A second scaling limiter then protects positivity in parameter space. For the density and energy traces, let $\bar{v}_{\mathbf{j}}$ denote the parameter-cell average and $v^{\min}_{\mathbf{j}}$ $= \min_{\mathbf{q}} v(\mathbf{y}_{\mathbf{j},\mathbf{q}})$ the smallest of the four Gauss-node values. Computing $\theta_{\mathbf{j}}$ from these quantities exactly as in Eqn.~\ref{eqn:pos_lim_phys}, each Gauss-node value is replaced by
\begin{equation}
v(\mathbf{y}_{\mathbf{j},\mathbf{q}}) \leftarrow \bar{v}_{\mathbf{j}} + \theta_{\mathbf{j}}\big(v(\mathbf{y}_{\mathbf{j},\mathbf{q}}) - \bar{v}_{\mathbf{j}}\big).
\label{eqn:pos_lim_stoch}
\end{equation}

\subsection{Numerical flux}
\label{sec:num_flux}
The numerical flux is evaluated pointwise at every interface and Gauss node from the left and right states $\mathbf{U}_L$ and $\mathbf{U}_R$, with $\mathbf{F}_{L,R} = \mathbf{F}(\mathbf{U}_{L,R})$. 


Our default numerical flux is the local Lax--Friedrichs (Rusanov) flux \cite{rusanov1962,kuzmin2023property}:
\begin{equation}
\begin{split}
& a_{LR} = \max\{|u_L| + c_L,\; |u_R| + c_R\}, \\
& \hat{\mathbf{F}}^{\textrm{LxF}}(\mathbf{U}_L, \mathbf{U}_R) = \frac{1}{2}\big(\mathbf{F}_L + \mathbf{F}_R\big) - \frac{a_{LR}}{2}\big(\mathbf{U}_R - \mathbf{U}_L\big).
\end{split}
\label{eqn:LxF_flux}
\end{equation}
Since the Rusanov flux does not attempt to resolve the detailed wave pattern of the Riemann problem, it is somewhat diffusive but very robust in the presence of strong shocks and large parameter-induced variations.

\subsection{Boundary and initial conditions}
\label{sec:ic_bc}
Boundary conditions are imposed through \emph{ghost cells} which supply the extra neighbors needed by the three-point WENO3 stencil. 

Reflective conditions are imposed by filling each ghost cell with its mirror image across the wall:
\begin{equation}
\bar{\mathbf{U}}_{1-\ell,\mathbf{j}} = \mathbf{Q}\,\bar{\mathbf{U}}_{\ell,\mathbf{j}}, \quad
\bar{\mathbf{U}}_{N+\ell,\mathbf{j}} = \mathbf{Q}\,\bar{\mathbf{U}}_{N+1-\ell,\mathbf{j}}, \quad
\ell = 1, 2, \quad
\mathbf{Q} = \operatorname{diag}(1,-1,1),
\label{eqn:gc_bcs}
\end{equation}
which flips the sign of the momentum and imposes zero normal velocity at $x_L(\mathbf{y})$ and $x_R(\mathbf{y})$ for every realization. 

Our test case is the Sod shock tube \cite{sod1978survey}, a standard benchmark for compressible flow solvers. A closed tube is divided by a membrane into a high-pressure, high-density region on the left and a low-pressure, low-density region on the right. At $t = 0$ the membrane is removed, and three waves emerge: a shock wave moving right into the low-pressure gas, a \emph{contact discontinuity} following behind it, and a \emph{rarefaction} moving left. In physical coordinates, the initial primitive state is
\begin{equation}
(\rho, u, p)\big|_{t=0} = \begin{cases}
(1,\; 0,\; 1), & x < \tfrac{1}{2}, \\
(0.125,\; 0,\; 0.1), & x \geq \tfrac{1}{2},
\end{cases}
\qquad \rho e = \frac{p}{\gamma - 1},
\label{eqn:sod_ic}
\end{equation}
where the total energy reduces to internal energy since the gas is initially at rest. Letting $\mathbf{U}_0(x)$ denote the corresponding conserved state, the initial cell averages are computed by evaluating $\mathbf{U}_0$ at the cell center in $\xi$ and averaging over the Gauss nodes in $\mathbf{y}$:
\begin{equation}
\bar{\mathbf{U}}_{i,\mathbf{j}}(0) = \sum_{\mathbf{q}} w_{\mathbf{q}}\; \mathbf{U}_0\big(x(\xi_i; \mathbf{y}_{\mathbf{j},\mathbf{q}})\big).
\label{eqn:ref_ICs}
\end{equation}
Since the membrane is fixed in physical space, its position in reference space varies with $\mathbf{y}$, as shown by Eqn.~\ref{eqn:xi_disc}. The geometric parameters thus affect the solution in two ways: through the initial data and through the Jacobian factor $\alpha(\mathbf{y})$.

\subsection{Time integration and step-size selection}
\label{sec:time_int}
Eqn.~\ref{eqn:fv_fluxdiff} is advanced in time by the three-stage strong-stability-preserving Runge--Kutta method (SSP-RK3) \cite{shu1988ssp,gottlieb2001ssp,gottlieb2009ssp}, whose stages are given in TT form in Eqn.~\ref{eqn:tt_ssprk3}. Since the parameter reconstruction couples neighboring parameter cells, all realizations advance with a common step size, set by the CFL condition for the fastest characteristic speed $|u|+c$ of the Euler equations \cite{evans2010pde,leveque1992numerical}, multiplied by $\alpha$ in reference coordinates:
\begin{align}
a_{\max} &:= \max_{i,\,\mathbf{j},\,\mathbf{q}}\big(|u| + c\big), \qquad \alpha_{\max} := \max_{\mathbf{j}}\, \alpha(\mathbf{y}_{\mathbf{j}}), \label{eqn:max_signalspeed_ref}\\
\Delta t^{n} &= \min\Big\{\textrm{CFL}\,\frac{\Delta\xi}{\alpha_{\max}\,a_{\max}},\; T - t^{n}\Big\}, \label{eqn:CFL_ss_ref}
\end{align}
where $a_{\max}$ is evaluated from Gauss-node values obtained by applying the parameter-space reconstruction of Section~\ref{sec:weno_stoch} to the cell averages $\bar{\mathbf{U}}_{i,\mathbf{j}}$. $\textrm{CFL} \in (0,1)$ is prescribed by the user. 
%
%
\section{Tensor-Train Representation of the Stochastic Finite Volume Scheme}
\label{sec:tt}
In this section we review the tensor-train (TT) format to the extent needed for our method, and we explain how each stage of the operator $\mathcal{L}$ and of the SSP-RK3 update is carried out in TT arithmetic. 
We explain how the operator $\mathcal{L}$ and the SSP-RK3 update are carried out in TT arithmetic for a single train that covers a contiguous range of reference cells. The partition of the reference mesh into such ranges is the subject of Section~\ref{sec:blockwise}.

\subsection{The tensor-train format and its operations}
\label{sec:tt_format}\label{sec:tt_ops}
A \emph{tensor train} \cite{oseledets2011tt} represents an array $\mathsf{A} \in \mathbb{R}^{n_1 \times \cdots \times n_\kappa}$ as a product of matrices,
\begin{equation}
\mathsf{A}(k_1, \dots, k_\kappa) = \mathbf{C}_1(k_1)\, \mathbf{C}_2(k_2) \cdots \mathbf{C}_\kappa(k_\kappa), \quad \mathbf{C}_\ell(k_\ell) \in \mathbb{R}^{r_{\ell-1} \times r_\ell}, \quad r_0 = r_\kappa = 1,
\label{eqn:tt_def}
\end{equation}
with \emph{cores} $\mathbf{C}_\ell$ of size $r_{\ell-1} \times n_\ell \times r_\ell$ and \emph{TT ranks} $r_1, \dots, r_{\kappa-1}$. We use the following standard operations \cite{oseledets2011tt}. A mode-$\ell$ product $\mathsf{A} \times_\ell \mathbf{P}$ with a matrix $\mathbf{P}$ changes only the $\ell$-th core, $\mathbf{C}_\ell(k) \leftarrow \sum_{k'} \mathbf{P}(k,k')\, \mathbf{C}_\ell(k')$, and is exact and rank-preserving. When $\mathbf{P}$ has a single unit entry per row, it selects, repeats, reverses, or restricts indices; we write $\mathsf{A}|_{a:b}$ for the restriction of the first mode to positions $a, \dots, b$ (reversed if $a > b$). Sums, the Hadamard product $\odot$, and the concatenation $\operatorname{cat}_1$ of trains along the first mode are also exact, but the ranks add for sums and concatenations and multiply for Hadamard products. Each is therefore followed by TT-SVD \emph{rounding} $\mathcal{T}$, which truncates such that $\|\mathsf{A} - \mathcal{T}(\mathsf{A})\|_F \le \varepsilon_{\mathrm{tt}} \|\mathsf{A}\|_F$ unless a rank cap $r_{\max}$ is reached.

The limiters, the WENO weights, and the numerical fluxes are nonlinear and have no exact TT counterpart. For a function $f: \mathbb{R}^p \to \mathbb{R}$ and trains $\mathsf{A}_1, \dots, \mathsf{A}_p$ of identical shape, we approximate the array
\begin{equation}
\mathsf{B}(k_1, \dots, k_\kappa) = f\big(\mathsf{A}_1(k_1, \dots, k_\kappa), \dots, \mathsf{A}_p(k_1, \dots, k_\kappa)\big)
\label{eqn:tt_pointwise}
\end{equation}
by cross approximation \cite{oseledets2010cross}, which samples entries of $\mathsf{B}$ at index sets chosen by maximal-volume pivoting \cite{goreinov1997pseudo,mikhalev2018maxvol} and adapts the ranks by DMRG-type sweeps \cite{savostyanov2011dmrg}. Each sample needs only the corresponding entries of the inputs, so the full array is never assembled, and since $f$ acts on entries rather than cores, the result does not depend on the non-unique choice of cores. We write $\mathcal{X}[f;\, \mathsf{A}_1, \dots, \mathsf{A}_p]$ for the rounded result. Arrays given directly by a formula in their indices, such as the initial condition and the geometric factor, are approximated in the same way \cite{dolgov2020parallel}.

All TT operations are performed with the BoBa library \cite{bobalibrary}. Cross approximation is started from rank two, uses a rank increment of one, and is stopped at relative tolerance $\varepsilon_{\mathrm{tt}}$. Rounding uses $\varepsilon_{\mathrm{tt}} = 10^{-4}$ and $r_{\max} = 128$. Every nonlinear step of the scheme is a cross approximation of the form of Eqn.~\ref{eqn:tt_pointwise}; no nonlinear operation is applied to individual cores.

\subsection{Stochastic finite volume unknowns as tensor trains}
\label{sec:tt_sfv}
Let $d$ be the number of uncertain parameters ($d = 2$ in Section~\ref{sec:FV_scheme}, with directions $1 \equiv L$ and $2 \equiv R$), and let $P = M + 4$ be the number of parameter cells per direction, including padding. With cell indices $\mathbf{j} = (j_1, \dots, j_d)$ and Gauss-node indices $\mathbf{q} = (q_1, \dots, q_d)$, $q_m \in \{1, 2\}$, the Gauss nodes and weights in Eqns.~\ref{eqn:sem_dis_geounc} and~\ref{eqn:ref_ICs} are the tensor product of the two-point rule, $(\mathbf{y}_{\mathbf{j},\mathbf{q}})_m = y_{m,j_m} + h_m\,\zeta_{q_m}$ and $w_{\mathbf{q}} = 2^{-d}$, with $\zeta_{1,2}$ as in Section~\ref{sec:weno_stoch}.

For a contiguous range of reference cells $\{i_0, \dots, i_0 + n - 1\}$, the values of the $\nu$-th conserved variable form the array $\bar{\mathsf{W}}_\nu \in \mathbb{R}^{n \times P \times \cdots \times P}$, which we store as a train with one physical and $d$ parameter cores,
\begin{equation}
\bar{\mathsf{W}}_\nu(\iota, j_1, \dots, j_d) := \big(\bar{\mathbf{U}}_{i_0 + \iota - 1,\, \mathbf{j}}\big)_\nu = \mathbf{C}_\xi(\iota)\, \mathbf{C}_1(j_1) \cdots \mathbf{C}_d(j_d), \qquad \iota = 1, \dots, n.
\label{eqn:sfv_tt}
\end{equation}
Its ranks satisfy $r_m \le \min\{n P^{m-1},\, P^{d-m+1}\}$, where $r_m$ links mode $m$ to mode $m+1$; in particular, the leading rank satisfies $r_1 \le \min\{n, P^d\}$, and for $d = 2$ the link between the two parameter cores satisfies $r_2 \le P$. Each conserved variable is a separate train, and we write $\bar{\mathsf{U}} = (\bar{\mathsf{W}}_1, \bar{\mathsf{W}}_2, \bar{\mathsf{W}}_3)$. Since the reference cell is the first mode, shifts, restrictions, reversals, and concatenations in $\xi$ act on $\mathbf{C}_\xi$ alone.

Several steps of the scheme act on point values at the Gauss nodes. We store these as trains of the same form on a \emph{node grid}, in which the cell index $j_m$ and node index $q_m$ are combined into $\hat{\jmath}_m = 2(j_m - 1) + q_m \in \{1, \dots, 2P\}$. Copying each cell value to both nodes, selecting node $q$, and Gauss averaging are mode products with
\begin{equation}
\mathbf{E} = \mathbf{I}_{P} \otimes \mathbf{1}_2, \qquad \mathbf{P}_q = \mathbf{I}_{P} \otimes \mathbf{b}_q^{\top}, \qquad \tfrac{1}{2}\mathbf{E}^{\top},
\label{eqn:node_maps}
\end{equation}
respectively, where $\mathbf{1}_2 = (1,1)^{\top}$ and $\mathbf{b}_q$ is the $q$-th canonical basis vector of $\mathbb{R}^2$. They are exact, rank-preserving, and act on one parameter mode at a time, so different parameter modes may lie on different grids.

The geometric factor is computed once by cross approximation as a $d$-mode train $\mathsf{A}_\alpha$ on the node grid with entries $\alpha(\mathbf{y}_{\mathbf{j},\mathbf{q}})$. Prepending a first core of size $1 \times n' \times 1$ with all entries equal to one gives the train $\mathsf{A}_\alpha^{[n']}$ that is constant over $n'$ reference interfaces and has the same ranks. The initial trains are cross approximations of the right-hand side of Eqn.~\ref{eqn:ref_ICs} at $\xi_{i_0 + \iota - 1}$.

\subsection{Evaluating the semi-discrete operator in tensor-train format}
\label{sec:tt_operator}
The three-point stencil of Section~\ref{sec:weno_phys} reaches two cells beyond the outermost interfaces of the range. The operator $\mathcal{L}$ of Eqn.~\ref{eqn:fv_fluxdiff} is therefore evaluated on trains $\bar{\mathsf{U}}$ whose first mode covers the extended range $\{i_0 - 2, \dots, i_0 + n + 1\}$, assembled as described in Section~\ref{sec:block_rhs}. Unless stated otherwise, each step is applied to each component separately, and we suppress the index $\nu$.

\textbf{\emph{Physical reconstruction.}} Let $\mathbf{S}_{\pm}$ be the shift matrices with unit entries $\mathbf{S}_-(\iota, \max\{\iota - 1, 1\})$ and $\mathbf{S}_+(\iota, \min\{\iota + 1, n+4\})$, so that $\bar{\mathsf{W}} \times_1 \mathbf{S}_{\pm}$ holds $\bar{w}_{i\pm1}$ at the position of cell $i$; the faces affected by the truncated shifts at the two ends are discarded. Let $\phi^-(z_-, z_0, z_+)$ denote the right-face value $w^-_{i+\frac{1}{2}}$ of Eqn.~\ref{eqn:right_recon} for the data $(\bar{w}_{i-1}, \bar{w}_i, \bar{w}_{i+1}) = (z_-, z_0, z_+)$, and let $\phi^+(z_-, z_0, z_+)$ denote the corresponding left-face value $w^+_{i-\frac{1}{2}}$. Both are pointwise functions of three values, so the face values of all cells and all parameter cells are obtained by two cross approximations:
\begin{equation}
\mathsf{W}^{\pm} = \mathcal{X}\big[\phi^{\pm};\; \bar{\mathsf{W}} \times_1 \mathbf{S}_-,\; \bar{\mathsf{W}},\; \bar{\mathsf{W}} \times_1 \mathbf{S}_+\big].
\label{eqn:tt_weno_phys}
\end{equation}
The density and energy, are then limited as in Eqn.~\ref{eqn:pos_lim_phys} by
\begin{equation}
\Theta = \mathcal{X}\big[\vartheta\big(\bar{w}, \min\{w^+, w^-\}\big);\; \bar{\mathsf{W}},\, \mathsf{W}^+,\, \mathsf{W}^-\big], \qquad \mathsf{W}^{\pm} \leftarrow \mathcal{T}\Big(\bar{\mathsf{W}} + \Theta \odot \big(\mathsf{W}^{\pm} - \bar{\mathsf{W}}\big)\Big),
\label{eqn:tt_pos_phys}
\end{equation}
where $\vartheta$ is the limiter factor including the clipping to $[0,1]$. By Section~\ref{sec:weno_phys}, the left and right states at the $n+1$ interfaces of the range are the restrictions $\bar{\mathsf{V}}^{L} = \mathsf{W}^-|_{2:n+2}$ and $\bar{\mathsf{V}}^{R} = \mathsf{W}^+|_{3:n+3}$, whose $\lambda$-th entry belongs to $\xi_{i_0 + \lambda - \frac{3}{2}}$ and is computed from a complete stencil.

\textbf{\emph{Parameter reconstruction.}}
Let $\bar{\mathsf{V}}$ denote any one of the six interface trains (three
components, two sides). The reconstruction of Section~\ref{sec:weno_stoch} is applied one parameter direction at
a time. Each sweep replaces the cell index $j$ of one mode by the node index
$\hat{\jmath} = 2(j-1) + q$, $q \in \{1,2\}$. Let $\mathbf{E}_0, \mathbf{E}_\pm \in \mathbb{R}^{2P \times P}$ map
cell averages to Gauss nodes. Each row $\hat{\jmath} = 2(j-1)+q$ contains a
single unit entry, so node $\hat{\jmath}$ copies one cell average:
$\bar v_j$ for $\mathbf{E}_0$, and $\bar v_{j\pm1}$ for $\mathbf{E}_\pm$.
In the boundary cells $j \in \{1, P\}$, which have no neighbor,
$\mathbf{E}_\pm$ copy $\bar v_j$ instead. Hence the trains
$\mathsf{A} \times_{m+1} \mathbf{E}_-$, $\mathsf{A} \times_{m+1} \mathbf{E}_0$
and $\mathsf{A} \times_{m+1} \mathbf{E}_+$ hold, at every node in direction
$y_m$, the stencil averages $\bar v_{j-1}$, $\bar v_j$ and $\bar v_{j+1}$.
Each product only copies slices of the $(m+1)$-th core, so it is exact and
leaves the TT ranks unchanged. At each node in direction $y_m$, they hold the averages
$\bar{v}_{j-1}$, $\bar{v}_j$ and $\bar{v}_{j+1}$ of the cell $j$ containing
that node. The position dependence of the ideal weights in
Eqn.~\ref{eqn:weno_stoch_lin} is supplied by a rank-one train $\mathsf{Z}_m$
with the same shape as these products. Its entry is $\zeta_q$ at node index
$\hat{\jmath}_m$, and it is constant along all other modes. With
$\psi(z_-, z_0, z_+, \zeta)$ the right-hand side of Eqn.~\ref{eqn:weno_stoch} for
$(\bar{v}_{j-1}, \bar{v}_j, \bar{v}_{j+1}) = (z_-, z_0, z_+)$, the sweeps start from
$\mathsf{V}_{(0)} = \bar{\mathsf{V}}$ and read
\begin{equation}
\mathsf{V}_{(m)} = \mathcal{X}\big[\psi;\;
  \mathsf{V}_{(m-1)} \times_{m+1} \mathbf{E}_-,\;
  \mathsf{V}_{(m-1)} \times_{m+1} \mathbf{E}_0,\;
  \mathsf{V}_{(m-1)} \times_{m+1} \mathbf{E}_+,\;
  \mathsf{Z}_m\big],
\label{eqn:tt_weno_stoch}
\end{equation}
for $m = 1 \dots d$. After sweep $m$, the train holds node values in $y_1, \dots, y_m$ and cell
averages in $y_{m+1}, \dots, y_d$. In the final train
$\mathsf{V} = \mathsf{V}_{(d)}$, the entry
$(\lambda, \hat{\jmath}_1, \dots, \hat{\jmath}_d)$ approximates the point value
of the interface state at $\mathbf{y}_{\mathbf{j},\mathbf{q}}$, where
$\mathbf{j} = (j_1, \dots, j_d)$ and $\mathbf{q} = (q_1, \dots, q_d)$. At the
outermost padding cells $j \in \{1, P\}$, the three shifted inputs coincide,
and $\psi$ returns the cell average. This is the constant fallback of
Section~\ref{sec:weno_stoch}.

\textbf{\emph{Flux, quadrature, and flux difference.}} Each component $\hat{F}_\nu$ of the numerical flux, is a pointwise function of the six node-valued states and is obtained by its own cross approximation. The geometric factor is applied at the nodes by an exact Hadamard product, and the Gauss average with $w_{\mathbf{q}} = 2^{-d}$ is an exact mode product:

\begin{equation}
\mathsf{F}_\nu = \mathcal{X}\big[\hat{F}_\nu;\; \mathsf{V}^{L}_1, \mathsf{V}^{L}_2, \mathsf{V}^{L}_3, \mathsf{V}^{R}_1, \mathsf{V}^{R}_2, \mathsf{V}^{R}_3\big], \qquad
\mathsf{G}_\nu = \mathcal{T}\big(\mathsf{F}_\nu \odot \mathsf{A}_\alpha^{[n+1]}\big) \times_2 \tfrac{1}{2}\mathbf{E}^{\top} \cdots \times_{d+1} \tfrac{1}{2}\mathbf{E}^{\top}.
\label{eqn:tt_flux}
\end{equation}
Thus $\alpha$ stays inside the parameter integral, as required in Section~\ref{sec:sfv_grids}, and the $\lambda$-th entry of $\mathsf{G}_\nu$ approximates the $\nu$-th component of $\mathbf{G}_{i_0 + \lambda - \frac{3}{2},\, \mathbf{j}}$. Finally, Eqn.~\ref{eqn:fv_fluxdiff} is a difference of two restrictions,
\begin{equation}
\mathsf{L}_\nu(\bar{\mathsf{U}}) = \mathcal{T}_{\mathrm{c}}\Big(-\frac{1}{\Delta\xi}\big(\mathsf{G}_\nu|_{2:n+1} - \mathsf{G}_\nu|_{1:n}\big)\Big),
\label{eqn:tt_fluxdiff}
\end{equation}
which approximates the $\nu$-th component of $\mathcal{L}(\bar{\mathbf{U}})$ on the cells of the range, and $\mathcal{T}_{\mathrm{c}}$ is the mean-conservative rounding described in Section~\ref{sec:tt_cons_round}.

\subsection{Time integration in tensor-train format}
\label{sec:tt_time}
Each stage of SSP-RK3 \cite{shu1988ssp} is formed as an exact sum and rounded once $\mathcal{T}_{\mathrm{c}}$:
\begin{equation}
\begin{split}
\bar{\mathsf{W}}^{(1)}_\nu &= \mathcal{T}_{\mathrm{c}}\Big(\bar{\mathsf{W}}^{n}_\nu + \Delta t\,\mathsf{L}_\nu(\bar{\mathsf{U}}^{n})\Big), \\
\bar{\mathsf{W}}^{(2)}_\nu &= \mathcal{T}_{\mathrm{c}}\Big(\tfrac{3}{4}\bar{\mathsf{W}}^{n}_\nu + \tfrac{1}{4}\bar{\mathsf{W}}^{(1)}_\nu + \tfrac{1}{4}\Delta t\,\mathsf{L}_\nu(\bar{\mathsf{U}}^{(1)})\Big), \\
\bar{\mathsf{W}}^{n+1}_\nu &= \mathcal{T}_{\mathrm{c}}\Big(\tfrac{1}{3}\bar{\mathsf{W}}^{n}_\nu + \tfrac{2}{3}\bar{\mathsf{W}}^{(2)}_\nu + \tfrac{2}{3}\Delta t\,\mathsf{L}_\nu(\bar{\mathsf{U}}^{(2)})\Big).
\end{split}
\label{eqn:tt_ssprk3}
\end{equation}
For $a_{\max}$ in Eqn.~\ref{eqn:max_signalspeed_ref}, Eqn.~\ref{eqn:tt_weno_stoch} is applied to $\bar{\mathsf{W}}_\nu$ without limiting, and the node values at the interior parameter cells are evaluated from the cores by Eqn.~\ref{eqn:tt_def}. The factor $\alpha_{\max}$ is computed directly at the parameter-cell centers.

\subsection{\texorpdfstring{Mean-conservative rounding}{Mean-conservative rounding}}
\label{sec:tt_cons_round}
The scheme of Section~\ref{sec:FV_scheme} is conservative for every realization, but every rounding $\mathcal{T}$ perturbs the cell averages and hence the totals of the conserved variables, and these perturbations accumulate over many steps. We therefore round the solution with a modified rounding routine that keeps the parameter mean of each total up to round-off. The idea is to round as usual and
then add back the lost total, using the smallest correction that does not
increase the TT ranks.

Let $\boldsymbol{\pi}$ $\in \mathbb{R}^{P}$ hold the averaging weights in one parameter direction: $1/M$ at the interior parameter cells and zero at the padding cells. For a train $\mathsf{A}$ of the form of Eqn.~\ref{eqn:sfv_tt},we
define the weighted total by summing over the spatial cells and averaging over
the interior parameter cells,
\begin{equation}
\begin{split}
\langle \mathsf{A} \rangle &:= \mathsf{A} \times_1 \mathbf{1}_n^{\top}
  \times_2 \boldsymbol{\pi}^{\top} \cdots \times_{d+1} \boldsymbol{\pi}^{\top}
  = \sum_{\iota=1}^{n} \mathbf{C}_\xi(\iota)\, \boldsymbol{\upsilon}, \\
\boldsymbol{\upsilon} &:= \sum_{\mathbf{j}} \pi_{j_1} \cdots \pi_{j_d}\,
  \mathbf{C}_1(j_1) \cdots \mathbf{C}_d(j_d) \in \mathbb{R}^{r_1}.
\end{split}
\label{eqn:mean_total}
\end{equation}
The mean-conservative rounding $\mathcal{T}_{\mathrm{c}}(\mathsf{A})$ forms $\mathcal{T}(\mathsf{A})$ and the deficit $\Delta_{\mathrm{c}} = \langle \mathsf{A} \rangle - \langle \mathcal{T}(\mathsf{A}) \rangle$, right-orthogonalizes the parameter cores of $\mathcal{T}(\mathsf{A})$, and restores the deficit by a uniform shift of the physical core,
\begin{equation}
\mathbf{C}_\xi(\iota) \leftarrow \mathbf{C}_\xi(\iota) + \frac{\Delta_{\mathrm{c}}}{n\,\|\boldsymbol{\upsilon}\|^2}\, \boldsymbol{\upsilon}^{\top}, \qquad \iota = 1, \dots, n,
\label{eqn:cons_correction}
\end{equation}
with $\boldsymbol{\upsilon}$ computed from the orthogonalized cores. This is the smallest change of $\mathbf{C}_\xi$ that restores the total, and it adds no rank. Since the parameter cores are orthonormal, the train changes by $|\Delta_{\mathrm{c}}|/(\sqrt{n}\,\|\boldsymbol{\upsilon}\|)$. This is never smaller than $|\Delta_{\mathrm{c}}|/(\sqrt{n}\,\|\boldsymbol{\pi}\|^{d})$, the change caused by the
smallest correction overall, the rank-one correction, $\Delta_{\mathrm{c}}\,(\mathbf{1}_n \otimes \boldsymbol{\pi} \otimes \cdots \otimes \boldsymbol{\pi})/(n\,\|\boldsymbol{\pi}\|^{2d})$. If the shift is more than twice as large, which happens when $\boldsymbol{\pi}$ is poorly represented by the retained parameter cores, the rank-one correction is added instead and the result is rounded at round-off level. In either case the correction is of the size of the rounding error, and a deficit at round-off level is left uncorrected.

We use $\mathcal{T}_{\mathrm{c}}$ for every rounding of the solution trains and of the flux differences, i.e., in Eqns.~\ref{eqn:tt_fluxdiff} and~\ref{eqn:tt_ssprk3} and in the split and merge of Section~\ref{sec:split_merge}, and $\mathcal{T}$ inside the operator. The parameter cells are not conserved separately, and the total of each individual parameter cell drifts at the level of the TT tolerance.
%
%
\section{Adaptive Blockwise Tensor-Train Representation}
\label{sec:blockwise}
We now partition the reference mesh into contiguous ranges of cells, called \emph{blocks}, each stored as a train as in Section~\ref{sec:tt}, and explain how neighboring blocks are coupled and how the partition is adapted. The reference mesh, the parameter grid, and the discretization of Section~\ref{sec:FV_scheme} are unaffected by the partition.

\subsection{Block partition}
\label{sec:block_partition}
We partition the reference cells $\{1, \dots, N\}$ into $B$ contiguous index sets
\begin{equation}
\mathcal{I}_b = \{i_b^-, \dots, i_b^+\}, \qquad b = 1, \dots, B, \qquad i_1^- = 1, \quad i_B^+ = N, \quad i_{b+1}^- = i_b^+ + 1,
\label{eqn:block_partition}
\end{equation}
with $n_b = |\mathcal{I}_b| = i_b^+ - i_b^- + 1$ cells in block $b$. Each block stores a triple $\bar{\mathsf{U}}^{(b)} = (\bar{\mathsf{W}}^{(b)}_1, \bar{\mathsf{W}}^{(b)}_2, \bar{\mathsf{W}}^{(b)}_3)$ of trains of the form of Eqn.~\ref{eqn:sfv_tt}, with $i_0 = i_b^-$, $n = n_b$, and their own cores and ranks, so that the leading rank satisfies $r_1 \le \min\{n_b, P^d\}$ (Section~\ref{sec:tt_sfv}). Concatenating the trains of all blocks along the first mode recovers the full SFV array. The global train corresponds to $B = 1$, and the per-cell trains to $n_b = 1$ for all $b$.

The adaptation routine of Section~\ref{sec:split_merge} is constrained by two block sizes. A block is \emph{split} only if both halves keep at least $n_{\min}$ cells. A block is \emph{merged} only it the result has a most $n_{\max}$ cells. The choice of $n_{\min} \ge 2$ guarantees that every block can supply the two halo cells require by its neighbors (see Section~\ref{sec:block_rhs}). The initial partition is uniform and may contain blocks larger than $n_{\max}$, which only affects merges. In the limit of a single initial block, the partition starts as the global TT and the adaptation builds blocks (see Section~\ref{sec:results_param6}). 


\subsection{Halo exchange and the block-local operator}
\label{sec:block_rhs}
The operator of Section~\ref{sec:tt_operator} requires, for block $b$, trains on the extended range $\{i_b^- - 2, \dots, i_b^+ + 2\}$. We call the two additional cells on each side the \emph{halo cells} of the block. Interior halo cells are copies of the edge cells of the neighboring blocks. At the walls, the halo cells are the ghost cells of Eqn.~\ref{eqn:gc_bcs}, built from the block's own edge cells. Let $Q_\nu$ $\in \{1, -1\}$ denote the $\nu$-th diagonal entry of the
reflection matrix $\mathbf{Q}$, so that
$Q_\nu$ $= -1$ for the velocity component and $Q_\nu$ $= 1$ otherwise. The left and right halo trains of block $b$ are,
\begin{equation}
\begin{split}
\mathsf{H}^{(b)}_{\nu,-} &=
\begin{cases}
\bar{\mathsf{W}}^{(b-1)}_\nu|_{n_{b-1}-1\,:\,n_{b-1}}, & b > 1, \\[2pt]
Q_\nu\, \bar{\mathsf{W}}^{(1)}_\nu|_{2:1}, & b = 1,
\end{cases} \\
\mathsf{H}^{(b)}_{\nu,+} &=
\begin{cases}
\bar{\mathsf{W}}^{(b+1)}_\nu|_{1:2}, & b < B, \\[2pt]
Q_\nu\, \bar{\mathsf{W}}^{(B)}_\nu|_{n_B\,:\,n_B-1}, & b = B.
\end{cases}
\label{eqn:halo}
\end{split}
\end{equation}
The reversed restrictions at the walls reproduce Eqn.~\ref{eqn:gc_bcs}. The extended trains are then assembled by concatenation along the first mode, followed by rounding:
\begin{equation}
\bar{\mathsf{W}}^{(b),\mathrm{ext}}_\nu = \mathcal{T}\Big(\operatorname{cat}_1\big(\mathsf{H}^{(b)}_{\nu,-},\; \bar{\mathsf{W}}^{(b)}_\nu,\; \mathsf{H}^{(b)}_{\nu,+}\big)\Big), \qquad \nu = 1, 2, 3.
\label{eqn:block_ext}
\end{equation}
All operations in Eqns.~\ref{eqn:halo} and~\ref{eqn:block_ext} other than the final rounding are exact. A halo train consists of two rows of the neighbor's $\mathbf{C}_\xi$ and the neighbor's parameter cores, which generally differ from those of block $b$; the concatenation places the two sets of parameter cores in block-diagonal position, and the rounding recompresses the extended train.

Applying the operator of Section~\ref{sec:tt_operator} to $\bar{\mathsf{U}}^{(b),\mathrm{ext}}$, with $i_0 = i_b^-$ and $n = n_b$, yields the block-local trains $\mathsf{L}^{(b)}_\nu$ on the cells $\mathcal{I}_b$. The flux at an interface shared by two blocks is computed within each of them from their respective extended trains, so the two values generally differ at the level of the TT tolerance. They are replaced by a single value as described in Section~\ref{sec:block_flux_share}.

\subsection{\texorpdfstring{Single-valued interface fluxes and conservation}{Single-valued interface fluxes and conservation}}
\label{sec:block_flux_share}
For the scheme to be conservative, the flux leaving one block through a shared interface must equal the flux entering its neighbor. We therefore evaluate the block-local operators in three steps. First, the interface fluxes $\mathsf{G}^{(b)}_\nu$ of Eqn.~\ref{eqn:tt_flux} are computed for every block. Second, at each interface shared by blocks $b$ and $b+1$, the last entry along the first mode of $\mathsf{G}^{(b)}_\nu$ and the first entry of $\mathsf{G}^{(b+1)}_\nu$ are both replaced by their average. The average is formed from the restrictions of the two original trains to a single interface, and it is rounded at round-off level. Each block then reassembles its flux train by concatenation of the two shared interfaces and the restriction to its own interior interfaces. At the two walls, the mass and energy components of the flux are set to zero, as the reflecting condition requires, while the momentum flux, which is the wall pressure, is kept. Third, the flux differences of Eqn.~\ref{eqn:tt_fluxdiff} are formed for every block. The second step uses only exact restrictions, sums, and concatenations, and exchanges one interface of each neighboring flux train.

With single-valued fluxes, the flux differences telescope across the blocks, so the weighted totals of Eqn.~\ref{eqn:mean_total} of their mass and energy components sum to zero over all blocks. Since $\mathcal{T}_{\mathrm{c}}$ keeps the weighted total of every train it rounds, both in the stage updates of Eqn.~\ref{eqn:tt_ssprk3} and in the split and merge of Section~\ref{sec:split_merge}, the parameter-mean total mass and energy are conserved up to round-off for any partition. The total momentum is not conserved even in exact arithmetic, since the walls exert a pressure on the gas.

\subsection{Time stepping over the partition}
\label{sec:block_time}
All blocks advance with the common step size $\Delta t$ of Eqn.~\ref{eqn:CFL_ss_ref}, where the maxima in Eqn.~\ref{eqn:max_signalspeed_ref} are taken over all blocks. Within a time step the partition is fixed, and the stages of Eqn.~\ref{eqn:tt_ssprk3} are carried out synchronously: at each stage, the operators $\mathsf{L}^{(b)}_\nu$ of all blocks are evaluated from the states of the preceding stage, including the halo cells of Eqn.~\ref{eqn:halo}, before any block is updated. The partition is adapted after the step has been completed, either after every step or after every fixed number of steps, which we call the adaptation interval, as described in the remainder of this section.

\subsection{Adaptation indicators}
\label{sec:block_indicators}
Adaptation is driven by three indicators per block, which are evaluated after each completed time step. Two of them are computed from cell-wise statistics of the solution: the mean density $\langle \rho \rangle_i$, the mean pressure $\langle p \rangle_i$, and the density variance $\sigma^2_{\rho,i}$ in reference cell $i$. These are taken over the Gauss nodes of the interior parameter cells, as defined in Section~\ref{sec:results_setup}, Eqn.~\ref{eqn:stats}. The third indicator depends only on the tensor train representations themselves.

The \emph{shock indicator} measures the largest relative jump of the mean density and mean pressure, and of the associated bands $\langle g \rangle_{i}\pm 1.96\sigma_{g,i}$, across the interfaces that bound or lie inside the block. Let
\begin{equation}
\mathcal{J}_b = \big\{\max\{1,\, i_b^- - 1\}, \dots, \min\{N - 1,\, i_b^+\}\big\}
\label{eqn:shock_ind_set}
\end{equation}
index the interfaces $\xi_{i+\frac{1}{2}}$ between two reference cells that touch $\mathcal{I}_b$. This set includes the two interfaces shared with the neighboring blocks, but no wall. The indicator is
\begin{equation}
\eta^{\mathrm{sh}}_b = \max_{g \in \{\rho, p \}}\; \max_{h \in \mathcal{H}_{g}}\; \max_{i \in \mathcal{J}_b} \frac{|h_{i+1} - h_i|}{\max\{|g_i|,\, 10^{-12}\}}, 
\label{eqn:shock_ind}
\end{equation}
where $\mathcal{H}_{g} = \{\langle g \rangle, \langle g \rangle - 1.96\sigma_{g}, \langle g \rangle + 1.96\sigma_{g} \}$. All jumps are normalized by the mean, since the bands can pass through zero. The jump of the mean alone will detect a wave only if the wave lies at similar reference positions in all realizations. If its reference position varies, averaging smears the jump of the mean over the range of positions, causing the mean to remain below any useful threshold even at a strong shock. In the test problem presented in Section~\ref{sec:results} relative jumps of the mean stay below about $0.11$. The bands, in contrast, jump sharply where the range of shock positions begins and ends, with relative jumps of up to about $0.7$. This allows the indicator to detect such an uncertain shock. The \emph{uncertainty indicator} responds instead to a spread of reference positions inside the band, which appears as variance in the density:
\begin{equation}
\eta^{\mathrm{unc}}_b = \bigg(\frac{1}{n_b} \sum_{i \in \mathcal{I}_b} \sigma^2_{\rho,i}\bigg)^{1/2}.
\label{eqn:unc_ind}
\end{equation}
Finally, the \emph{rank indicator} is the largest TT rank stored in the block,
\begin{equation}
\eta^{\mathrm{rk}}_b = \max_{\nu = 1, 2, 3}\, \max_{m = 1, \dots, d}\, r_m(\bar{\mathsf{W}}^{(b)}_\nu).
\label{eqn:rank_ind}
\end{equation}
For each indicator we prescribe a refinement threshold $\tau^{\uparrow}$ and a coarsening threshold $\tau^{\downarrow} < \tau^{\uparrow}$. The rank thresholds are set relative to $R = \min\{r_{\max}, P\}$, the bound on the link between the two parameter cores (Section~\ref{sec:tt_sfv}), which the block ranks attain in practice, so that they scale with the parameter resolution. By default,
\begin{equation}
\big(\tau^{\uparrow}_{\mathrm{sh}}, \tau^{\uparrow}_{\mathrm{unc}}, \tau^{\uparrow}_{\mathrm{rk}}\big) = (0.12,\, 0.04,\, [0.75R]), \quad
\big(\tau^{\downarrow}_{\mathrm{sh}}, \tau^{\downarrow}_{\mathrm{unc}}, \tau^{\downarrow}_{\mathrm{rk}}\big) = (0.03,\, 0.015,\, [0.4R]),
\label{eqn:ind_thresholds}
\end{equation}
where $[\cdot]$ denotes rounding to an integer value. 

\subsection{Splitting and merging}
\label{sec:split_merge}
 \textbf{\emph{Splitting.}} A block is marked for refinement if any of its indicators exceeds the corresponding refinement threshold,
\begin{equation}
\eta^{\mathrm{sh}}_b > \tau^{\uparrow}_{\mathrm{sh}} \quad \textrm{or} \quad \eta^{\mathrm{unc}}_b > \tau^{\uparrow}_{\mathrm{unc}} \quad \textrm{or} \quad \eta^{\mathrm{rk}}_b > \tau^{\uparrow}_{\mathrm{rk}},
\label{eqn:refine_crit}
\end{equation}
and if it is large enough to be halved without violating the lower size bound, i.e., $\lfloor n_b/2 \rfloor \ge n_{\min}$. A marked block is replaced by two adjacent blocks: the left one takes
the first $\lfloor n_b/2 \rfloor$ cells of $\mathcal{I}_b$ and the right
one the remaining $\lceil n_b/2 \rceil$, so that the right half receives
the additional cell when $n_b$ is odd. The new trains are the restrictions
\begin{equation}
\mathcal{T}_{\mathrm{c}}\big(\bar{\mathsf{W}}^{(b)}_\nu|_{1\,:\,\lfloor n_b/2 \rfloor}\big), \qquad \mathcal{T}_{\mathrm{c}}\big(\bar{\mathsf{W}}^{(b)}_\nu|_{\lfloor n_b/2 \rfloor + 1\,:\,n_b}\big), \qquad \nu = 1, 2, 3.
\label{eqn:split}
\end{equation}
Each restriction acts on $\mathbf{C}_\xi$ alone and is exact, but both halves would inherit the parameter cores of the parent, so each half is rounded. All marked blocks are split in a single pass, and each block is split at most once per adaptation.

\textbf{\emph{Merging.}} Coarsening is applied to pairs of adjacent blocks $b$ and $b+1$. A pair is a candidate for merging if every indicator lies below its coarsening threshold on both blocks,
\begin{equation}
\max\{\eta^{\mathrm{sh}}_b, \eta^{\mathrm{sh}}_{b+1}\} < \tau^{\downarrow}_{\mathrm{sh}}, \quad
\max\{\eta^{\mathrm{unc}}_b, \eta^{\mathrm{unc}}_{b+1}\} < \tau^{\downarrow}_{\mathrm{unc}}, \quad
\max\{\eta^{\mathrm{rk}}_b, \eta^{\mathrm{rk}}_{b+1}\} < \tau^{\downarrow}_{\mathrm{rk}},
\label{eqn:coarsen_crit}
\end{equation}
and if the merged block respects the upper size bound, $n_b + n_{b+1} \le n_{\max}$. For a candidate pair, the trains are concatenated exactly along the first mode and then rounded,
\begin{equation}
\mathsf{W}^{\mathrm{cat}}_\nu = \operatorname{cat}_1\big(\bar{\mathsf{W}}^{(b)}_\nu,\, \bar{\mathsf{W}}^{(b+1)}_\nu\big), \qquad \bar{\mathsf{W}}^{\mathrm{merge}}_\nu = \mathcal{T}_{\mathrm{c}}(\mathsf{W}^{\mathrm{cat}}_\nu), \qquad \nu = 1, 2, 3.
\label{eqn:merge}
\end{equation}
The concatenation places the parameter cores of the two blocks in block-diagonal position. Rounding then determines how much parametric structure the two blocks share. The merge is accepted only if rounding perturbs the density train by at most a prescribed tolerance,
\begin{equation}
\frac{\big\|\bar{\mathsf{W}}^{\mathrm{merge}}_1 - \mathsf{W}^{\mathrm{cat}}_1\big\|_F}{\|\mathsf{W}^{\mathrm{cat}}_1\|_F} \le \varepsilon_{\mathrm{merge}}, \qquad \varepsilon_{\mathrm{merge}} = 2 \times 10^{-3}.
\label{eqn:merge_test}
\end{equation}
Otherwise both blocks are kept. Since the rank cap $r_{\max}$ of Section~\ref{sec:tt_ops} applies in Eqn.~\ref{eqn:merge}, the test also rejects merges that the cap would truncate too strongly.

The pairs are examined in a sweep from left to right. After an accepted merge, the indicators of the merged block are evaluated immediately, Eqns.~\ref{eqn:shock_ind}--\ref{eqn:unc_ind} from the same cell-wise statistics and Eqn.~\ref{eqn:rank_ind} from $\bar{\mathsf{W}}^{\mathrm{merge}}_\nu$, and the merged block is compared with its right neighbor, so a smooth region can be coarsened over several blocks within one sweep. The indicators of all blocks are recomputed at the start of each sweep, from the cell-wise statistics computed once per adaptation. At most two sweeps are performed per adaptation, and the procedure stops early if a sweep accepts no merge. Splitting and merging only regroup the cells into trains and introduce no error beyond the roundings in Eqns.~\ref{eqn:split} and~\ref{eqn:merge}.

\subsection{Summary of one time step}
\label{sec:block_algorithm}
Algorithm~\ref{alg:block_step} summarizes one step of the adaptive blockwise solver: step-size selection, three synchronous SSP-RK3 stages on a fixed partition with halo exchange before and interface-flux exchange after every operator evaluation, and adaptation of the partition.

\begin{algorithm}
\caption{One time step of the adaptive blockwise TT solver.}
\label{alg:block_step}
\begin{algorithmic}[1]
\Require partition $\{\mathcal{I}_b\}_{b=1}^{B}$, trains $\bar{\mathsf{U}}^{(b),n}$ for $b = 1, \dots, B$, time $t^n$
\State Compute $a_{\max}$ over all blocks and $\alpha_{\max}$ as in Section~\ref{sec:tt_time}; set $\Delta t^n$ by Eqn.~\ref{eqn:CFL_ss_ref}
\For{each stage of Eqn.~\ref{eqn:tt_ssprk3}}
    \For{$b = 1, \dots, B$}
        \State Assemble $\bar{\mathsf{U}}^{(b),\mathrm{ext}}$ from the current stage states by Eqns.~\ref{eqn:halo}--\ref{eqn:block_ext}
        \State Evaluate the interface fluxes $\mathsf{G}^{(b)}_\nu$ of $\bar{\mathsf{U}}^{(b),\mathrm{ext}}$, $\nu = 1, 2, 3$, as in Section~\ref{sec:tt_operator}
    \EndFor
    \State Make the fluxes at shared interfaces single-valued and set the wall mass and energy fluxes to zero (Section~\ref{sec:block_flux_share})
    \State Form $\mathsf{L}^{(b)}_\nu$ for every block by Eqn.~\ref{eqn:tt_fluxdiff}
    \For{$b = 1, \dots, B$}
        \State Form and round the stage update of $\bar{\mathsf{W}}^{(b)}_\nu$, $\nu = 1, 2, 3$, by Eqn.~\ref{eqn:tt_ssprk3}
    \EndFor
\EndFor
\State Compute the cell-wise statistics and the indicators of Eqns.~\ref{eqn:shock_ind}--\ref{eqn:rank_ind}
\State Split every block that satisfies Eqn.~\ref{eqn:refine_crit} and $\lfloor n_b/2 \rfloor \ge n_{\min}$, by Eqn.~\ref{eqn:split}
\For{sweep $= 1, 2$}
    \State Recompute the indicators of all blocks from the same cell-wise statistics
    \State Sweep left to right over adjacent pairs; merge by Eqns.~\ref{eqn:coarsen_crit}--\ref{eqn:merge_test}, updating the indicators of each merged block
    \If{no merge was accepted}
        \State \textbf{break}
    \EndIf
\EndFor
\State $t^{n+1} = t^n + \Delta t^n$
\end{algorithmic}
\end{algorithm}

%
%
\section{\texorpdfstring{Numerical results}{Numerical results}}
\label{sec:results}
In this section we compare the adaptive blockwise TT with two classic layouts on the uncertain-geometry Sod problem of Section~\ref{sec:ic_bc}. We assess the accuracy of the statistics, the conservation of the mean, the ranks and block layout, and the cost on this problem, before extending the method to a case with eight uncertain parameters.

\subsection{\texorpdfstring{Setup}{Setup}}
\label{sec:results_setup}
\begin{figure}[htbp]
\centering
\includegraphics[width=\textwidth]{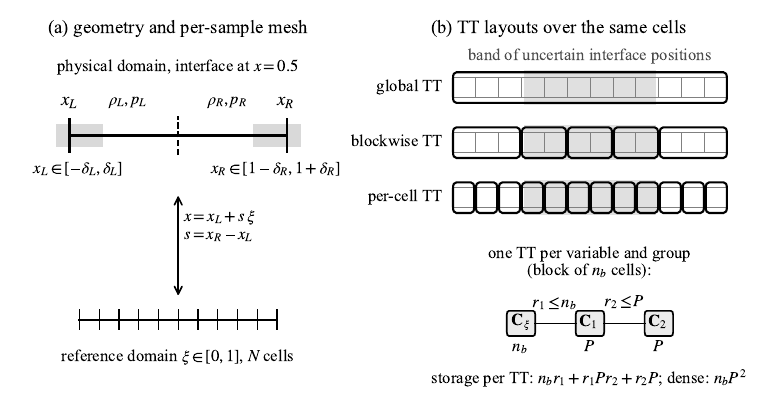}
\caption{Test problem and TT layouts. (a) Uncertain walls, fixed diaphragm, and the affine map to the reference interval. (b) Global, blockwise, and per-cell layouts over the same reference cells, with the band of uncertain diaphragm positions shaded, and the three cores of one train. The cores are $\mathbf{C}_\xi$, $\mathbf{C}_1$, $\mathbf{C}_2$ of Eqn.~\ref{eqn:sfv_tt}, with $P = M + 4$.}
\label{fig:schematic}
\end{figure}

Unless stated otherwise, the reference interval is divided into $N = 100$ cells, each parameter range into $M = 10$ cells, and $\delta_L = \delta_R = 0.24$. Since the diaphragm is fixed at $x = 1/2$ in physical space, its reference position $\xi_d(\mathbf{y})$ in Eqn.~\ref{eqn:xi_disc} covers a band of reference cells, which widens with $\delta_L$ and $\delta_R$ (Fig.~\ref{fig:schematic}(a)).

We compare three TT layouts (Fig.~\ref{fig:schematic}(b)). The \emph{global TT} stores each conserved variable in a single train ($B = 1$), and the \emph{per-cell TT} uses one train per conserved variable and reference cell ($n_b = 1$). The \emph{blockwise TT} of Section~\ref{sec:blockwise} starts from blocks of $10$ cells and, unless stated otherwise, uses $n_{\min} = 5$, the thresholds of Eqn.~\ref{eqn:ind_thresholds}, and is adapted after every fifth step. All runs use $\varepsilon_{tt}=10^{-4}$, except for the tolerance study of Fig.~\ref{fig:conservation}(c). All layouts use the same discretization and TT operations, including the mean-conservative rounding of Section~\ref{sec:tt_cons_round}, single-valued interface fluxes, and zero wall mass and energy fluxes, so all conserve the parameter-mean mass and energy. The rank cap $r_{\max}$ is reached only by the global TT on the finest mesh. To compare the layouts at identical time levels, the step size is not selected by Eqn.~\ref{eqn:CFL_ss_ref} but held fixed, using the upper bound on $\alpha_{\max}\,a_{\max}$ implied by $\alpha \le 1/s_{\min}$ (Eqn.~\ref{eqn:alpha_clamp}), and is shortened only to land on the output times. With $\delta_L = \delta_R = 0.24$, the clamp of Eqn.~\ref{eqn:alpha_clamp} acts only in the padding cells.

Statistics are computed from node values of the primitive variables. For each reference cell $i$, the parameter reconstruction of Section~\ref{sec:weno_stoch}, without the limiter, gives the conserved state at the Gauss nodes of every parameter cell, from which $\rho$, $u$, and $p$ are recovered. For $g \in \{\rho, u, p\}$, the mean and variance in cell $i$ are taken with equal weights over the four Gauss nodes $\mathbf{q}$ of each interior parameter cell $\mathbf{j}$,

\begin{equation}
\langle g \rangle_i = \frac{1}{4M^2} \sum_{\mathbf{j}} \sum_{\mathbf{q}} g_i(\mathbf{y}_{\mathbf{j},\mathbf{q}}), \qquad
\sigma^2_{g,i} = \frac{1}{4M^2} \sum_{\mathbf{j}} \sum_{\mathbf{q}} \big(g_i(\mathbf{y}_{\mathbf{j},\mathbf{q}}) - \langle g \rangle_i\big)^2.
\label{eqn:stats}
\end{equation}
These statistics also drive the adaptation indicators of Section~\ref{sec:block_indicators}. Alongside the profiles of $g \in \{\rho, u, p\}$ in Section~\ref{sec:results_accuracy} we report the central $95\%$ interval between the $2.5\%$ and $97.5\%$ percentiles of $g_i$, which are taken from the empirical distribution of the $4M^{2}$ Gauss-node values with equal weights. 

Statistics in physical coordinates are computed pointwise in $x$. Each Gauss node, $\mathbf{y}_{\mathbf{j},\mathbf{q}}$, defines a realization on its own physical mesh, whose cell interfaces lie at $x(\xi_{i\pm\frac12}; \mathbf{y}_{\mathbf{j},\mathbf{q}})$ and on which the realization is piecewise constant. To evaluate it at a physical point $x$, we first check whether $x$ lies in its domain $D_x(\mathbf{y}_{\mathbf{j},\mathbf{q}})$. If not, the node does not contribute. If so, we map $x$ back to
$\xi(x;\mathbf{y}_{\mathbf{j},\mathbf{q}})$ and take the value $g_i(\mathbf{y}_{\mathbf{j},\mathbf{q}})$ of the cell $K_i^{\xi}$ containing
this point. The sum of the weights of the contributing nodes approximates the coverage $\mu(\Gamma_x)$, i.e., the probability that a realization covers $x$. We renormalize these weights to sum to one and compute the conditional mean, variance, and percentiles at $x$ from the contributing nodes. We evaluate the statistics on a uniform grid with spacing $10^{-3}$ in $[-\delta_L, 1 + \delta_R]$, the union of all physical domains, and omit points with a coverage below $5\%$, where too few realizations contribute.

The dense reference solves Eqn.~\ref{eqn:cons_law_ref} separately at each point of a $400 \times 400$ tensor-product Gauss--Legendre grid in $\Gamma$, with the same physical discretization and time steps but without TT approximation, and computes the statistics with the Gauss--Legendre weights. Its difference from the TT statistics thus contains both the stochastic finite volume error and the TT truncation error. We measure it by the relative $L^1$ error over all reference cells, taken jointly over the three primitive variables,
\begin{equation}
e_{\mathrm{mean}} = \frac{\sum_{g} \sum_{i} \big|\langle g \rangle_i - \langle g \rangle^{\mathrm{ref}}_i\big|}{\sum_{g} \sum_{i} \big|\langle g \rangle^{\mathrm{ref}}_i\big|},
\label{eqn:err_rel}
\end{equation}
and define $e_{\mathrm{std}}$ for $\sigma_{g,i}$ analogously. As cost measures we report the peak storage, i.e., the largest number of entries held in the cores of all trains over all Runge--Kutta stages, compared with the $3N(M+4)^2$ entries of the dense SFV array, the wall time per Runge--Kutta step, and the number of trains. The scaling studies run a small, fixed number of steps from the initial condition and thus reflect the early phase, in which the band of uncertain diaphragm positions dominates the ranks.

\subsection{\texorpdfstring{Accuracy of the statistics}{Accuracy of the statistics}}
\label{sec:results_accuracy}

\begin{figure}[htbp]
\centering
\includegraphics[width=\textwidth]{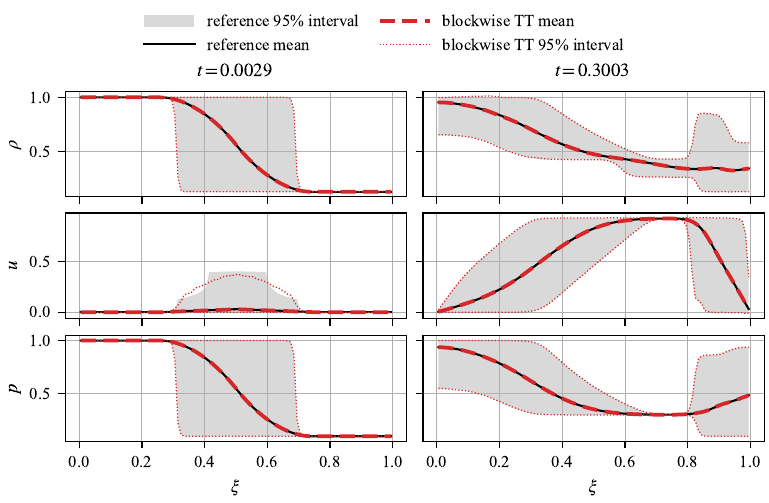} 
\caption{Mean and central $95\%$ interval, between the $2.5\%$ and $97.5\%$ percentiles over the
parameters, of $g = \rho$, $u$, and $p$ (rows) over the reference coordinate $\xi$ at the first and last output times. Reference: mean (solid black) and interval (gray band). Blockwise TT: mean (thick red dashed) and interval bounds (thin red dotted). $N = 100$, $M = 40$.}
\label{fig:profiles_mean}
\end{figure}
Fig.~\ref{fig:profiles_mean} compares the mean and the central 95\% interval of the blockwise TT with those of the reference at the first and last output time for a run with $M=40$. At these times, the blockwise mean differs from the reference mean by less than $0.5\%$ of the largest value of each variable, and the bounds of the two intervals agree to within about $2\%$. The exception is the upper bound of $u$ at the first output time. For $\xi$ between $0.4-0.6$ it differs from the upper bound of the reference by about $0.1\%$. At this time, the velocity in each realization is nonzero only in the narrow region (typically a few cells wide) around the location of the diaphragm. The diaphragm
position differs between the realizations and lies between $\xi \approx 0.4$
and $0.6$. Consider a fixed cell in this interval. The velocity in this cell is nonzero only for realizations whose diaphragm is at most a few cells away from it. Since the waves are narrow compared to the interval, only a small percentage of the realizations have a nonzero velocity in this cell. The upper
bound in this cell is thus determined by a small part of the parameter space, and its accuracy depends on how well the parameter discretization resolves
this part. The mean profiles are smooth, since averaging over the wave positions smears the jumps, and the interval moves with the waves, from the middle of the domain, where the diaphragm lies for some realizations but not for others, to the right wall, which the shock has reached in part of the ensemble by the last output time.

\begin{figure}[htbp]
\centering
\includegraphics[width=\textwidth]{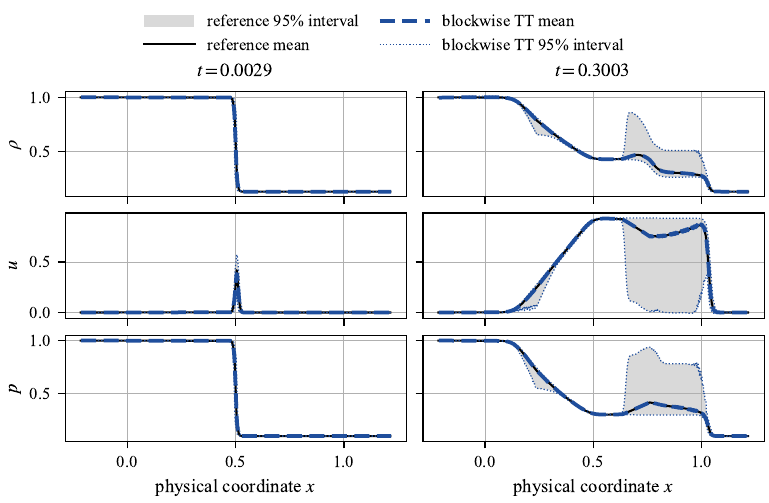} 
\caption{Mean and central $95\%$ interval of $g = \rho$, $u$, and $p$ as in Fig.~\ref{fig:profiles_mean}, but over the physical coordinate $x$, conditional on the realizations whose domain contains $x$ (Eqns.~\ref{eqn:cover_set}-\ref{eqn:phys_mean}).Reference: mean (solid black) and interval (gray band). Blockwise TT: mean (thick blue dashed) and interval bounds (thin blue dotted).
Same runs as in Fig. 2; $N = 100$, $M = 40$.}
\label{fig:profiles_mean_x}
\end{figure}

Fig.~\ref{fig:profiles_mean_x} shows the same statistics in physical coordinates (Section~\ref{sec:stat_map}). At the first three output times, the interval collapses onto the mean almost everywhere. This behavior is expected, as the diaphragm is fixed in physical space, and until a wave reaches a wall, the flow does not depend on the wall positions, so the solution at a fixed $x$ is nearly deterministic. The wide intervals in $\xi$ at these times (Fig.~\ref{fig:profiles_mean}) are thus created by the mapping alone. At the last output time, genuine uncertainty has appeared in two regions. Near the left wall (for $x$ between about $0.15$ and $0.35$), the rarefaction has been reflected from the wall for some realizations with large $x_{L}$. Near the right wall x (for $x$between about $0.65$ and $1.03$), the shock has been reflected in realizations with shorter domains but has not yet reached it in realizations with longer ones, so that the gas at a given $x$ is either still moving behind the incident shock, with $u\approx0.93$, or has been stopped and compressed by the reflected shock. Beyond $x\approx 1.03$, only the longer domains contribute, and
in all of them the gas is still at rest. Apart from the peaks at the first output time,
the blockwise means and intervals agree with the reference in
physical coordinates as closely as in reference coordinates (Fig.~\ref{fig:profiles_mean_x}).

\begin{figure}[htbp]
\centering
\includegraphics[width=\textwidth]{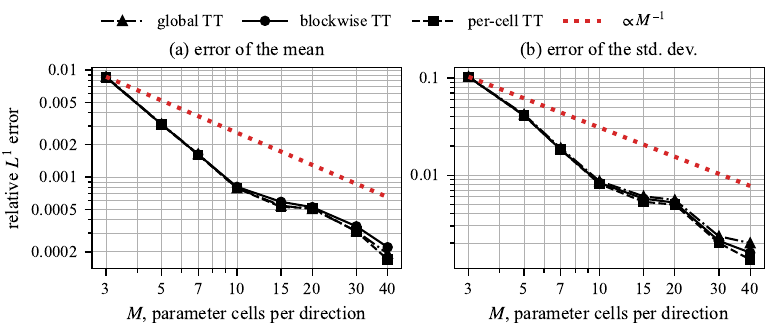}
\caption{Errors $e_{\mathrm{mean}}$ (left) and $e_{\mathrm{std}}$ (right) of Eqn.~\ref{eqn:err_rel} versus the number of parameter cells per direction, $M$, at the second output time, $t = 0.0507$, for the global, blockwise, and per-cell TT. Dotted: first order. $N = 100$.}
\label{fig:convergence}
\end{figure}

The errors $e_{\mathrm{mean}}$ and $e_{\mathrm{std}}$ decrease monotonically with $M$ at roughly first order or slightly faster (Fig.~\ref{fig:convergence}), as expected for a solution that is discontinuous in the parameters, and $e_{\mathrm{std}}$ is about an order of magnitude larger than $e_{\mathrm{mean}}$. At the first output time, the error of the three layouts coincide. At the second, shown in Fig.~\ref{fig:convergence}, they begin to diverge slightly after $M=10$. We attribute these
differences to the TT truncation, which in the blockwise and global layouts acts on several cells at once. They remain small compared with the discretization error.

\subsection{\texorpdfstring{Conservation}{Conservation}}
\label{sec:results_conservation}
\begin{figure}[htbp]
\centering
\includegraphics[width=\textwidth]{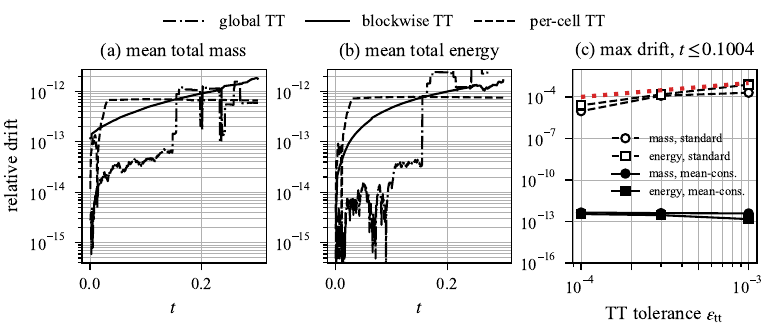}
\caption{Conservation of the parameter-mean totals. (a, b) Relative drift $|\mathcal{M}_\nu(t) - \mathcal{M}_\nu(0)|/|\mathcal{M}_\nu(0)|$ in time of the parameter-mean total mass and energy for the global, blockwise, and per-cell TT. (c) Largest relative drift of the parameter-mean total mass and energy over $t \in [0, 0.1004]$ versus the TT tolerance $\varepsilon_{\mathrm{tt}}$ for the blockwise TT with standard rounding and with the mean-conservative scheme of Sections~\ref{sec:tt_cons_round} and~\ref{sec:block_flux_share}. Dotted: $\varepsilon_{\mathrm{tt}}$; dash-dotted: machine precision. $N = 100$, $M = 10$.}
\label{fig:conservation}
\end{figure}

For all three layouts, the relative drift of the parameter-mean totals $\mathcal{M}_\nu(t) = \Delta\xi \sum_i \frac{1}{M^2} \sum_{\mathbf{j}} (\bar{\mathbf{U}}_{i,\mathbf{j}}(t))_\nu$ of mass and energy, over the interior parameter cells, stays within a few orders of magnitude of the machine precision over the whole run, many orders of magnitude below the TT tolerance (Fig.~\ref{fig:conservation}(a, b)). It stays below $4\times10^{-12}$ for every layout and variable. For the global TT it remains below about $2\times10^{-13}$ up to $t\approx0.15$ and then jumps to about $3\times10^{-12}$. For the per-cell TT the drift rises during the first steps, to about $7\times10^{-13}$, and then stays nearly constant. For the blockwise TT the drift grows slowly and steadily in
time, to $1.8\times10^{-12}$ at t= 0.3003, consistent with round-off accumulating over the many roundings and structural changes of the blocks, but it remains at round-off level, as guaranteed by the construction of Section~\ref{sec:block_flux_share}. In particular, the splits and merges of the blockwise TT do not compromise conservation.

To separate the effect of the mean-conservative scheme from that of the TT tolerance, Fig.~\ref{fig:conservation}(c) repeats the blockwise run up to $t = 0.1004$ for the tolerances $\varepsilon_{\mathrm{tt}}=10^{-3},3\times10^{-4}$ and $10^{-4}$, both with the mean-conservative scheme and with a standard scheme that uses the rounding $\mathcal{T}$ in place of $\mathcal{T}_{\mathrm{c}}$. With the standard scheme, the largest drift of the parameter-mean total mass and energy decreases roughly in proportion to $\varepsilon_{\mathrm{tt}}$, so conservation is an approximation controlled by the tolerance and can only be improved by tightening it, as observed for block-local TTs in radiation transport \cite{gorodetsky2026blocktt}. With the mean-conservative scheme, the drift stays at round-off level for every tolerance, including the loosest.

\subsection{\texorpdfstring{Ranks and block layout}{Ranks and block layout}}
\label{sec:results_rank}
\begin{figure}[htbp]
\centering
\includegraphics[width=\textwidth]{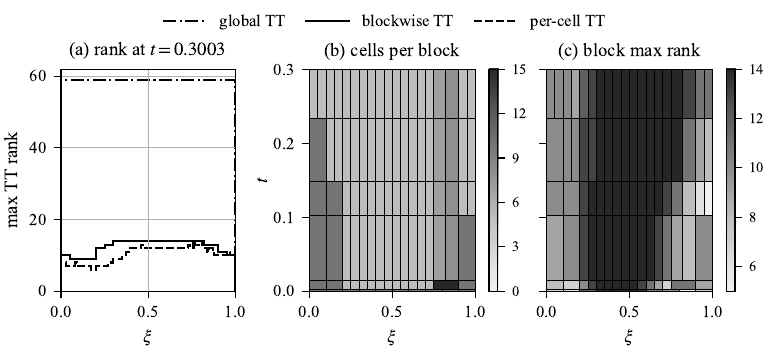}
\caption{(a) Largest TT rank over space at the last output time for the global, blockwise, and per-cell TT. (b) Number of cells per block and (c) largest rank per block of the blockwise TT in space and time. $N = 100$, $M = 10$.}
\label{fig:rank_layout}
\end{figure}

Fig.~\ref{fig:rank_layout}(a) shows the largest TT rank over space at the last output time, which for the blockwise and per-cell layouts is the largest rank of the trains covering each reference cell. The global rank, which the few cells containing the uncertain waves set for the whole domain, is several times larger than the local ranks everywhere. It is $74$ at $t = 0.05$, and at most $78$ during the run, whereas the blockwise and per-cell ranks never exceed $P = 14$. At $t = 0.05$, the blockwise and per-cell ranks are equal, at $P = 14$, in the middle of the domain, where the uncertain waves are, while towards the walls the blockwise rank is larger, $8$ against $4$ next to the left wall, since it also includes the coupling between the cells of a block through $r_1$. At the last output time, the blockwise rank is $14$ for $\xi$ between about $0.3$ and $0.85$ and $9$ to $10$ next to the walls, one to three above the per-cell rank. Panels (b) and (c) show how the partition responds to the waves. At the first adaptation, the blocks over the band of diaphragm positions, $0.2 \le \xi \le 0.75$, are split down to the minimum size of five cells. The partition changes again when the waves reach the walls: the block at the right wall is halved at $t \approx 0.10$, and the two blocks at the left wall at $t \approx 0.15$ and $t \approx 0.23$. By the last output time, all $19$ blocks have the minimum size. The block ranks reach $P = 14$ in the middle of the domain throughout and increase near the walls as the waves arrive.

\subsection{\texorpdfstring{Cost and storage}{Cost and storage}}
\label{sec:results_cost}
\begin{figure}[htbp]
\centering
\includegraphics[width=\textwidth]{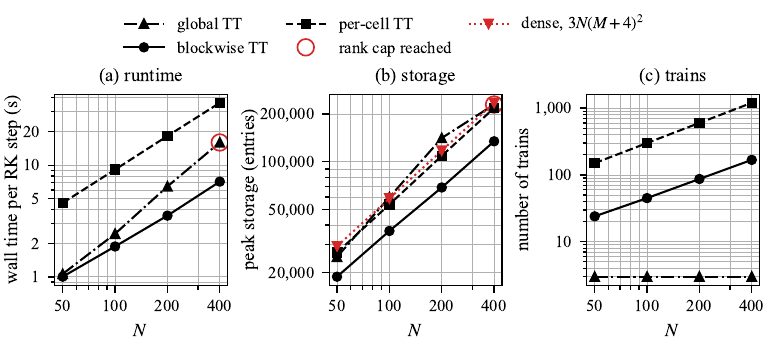} 
\caption{(a) Wall time per Runge--Kutta step, (b) peak storage, with the dense storage $3N(M+4)^2$, and (c) number of trains versus $N$ for the global, blockwise, and per-cell TT. Circles: rank cap $r_{\max}$ reached. $M = 10$, fixed number of steps.}
\label{fig:scaling_nx}
\end{figure}

With $M = 10$ fixed (Fig.~\ref{fig:scaling_nx}), the runtime per step grows roughly linearly in $N$ for the blockwise and per-cell TT and faster for the global TT, whose rank grows with $N$. The blockwise TT is the fastest layout at every $N$. On the finest mesh, the global TT is $2.2$ times slower than blockwise TT and the per-cell TT $5.1$ times slower. The blockwise TT is also the most compact layout at every $N$. At $N = 400$, its peak storage is $0.57$ times that of the dense array, $41\%$ below that of the global TT and $38\%$ below that of the per-cell TT, which both need nearly as much storage as the dense array. The blockwise TT uses about seven times fewer trains than the per-cell TT, whose number of trains is $3N$.

On the finest mesh, cross approximation accounts for most of the step time in all three layouts: $68\%$ for the global TT, $78\%$ for the blockwise TT, and $83\%$ for the per-cell TT. The global TT spends the largest share outside cross approximation, mostly on the rounding and bookkeeping of its high-rank trains. The number of TT operations per step is independent of $N$ for the global TT, and grows linearly with $N$ for the other two layouts. The blockwise TT performs about $40\%$ as many operations as the per-cell TT, but two orders of magnitude more than the global TT. Its runtime is nonetheless the smallest, because its operations act on trains with small ranks and short physical modes, whereas every operation of the global TT acts on a train whose rank is set by the uncertain waves. The per-cell TT is the slowest layout, since it performs $2.6$ times as many operations as the blockwise TT, which on average also take about twice as long.

When $M$ is varied at $N = 100$, the storage of all layouts grows more slowly than the dense storage, which grows like $(M+4)^2$. The blockwise TT needs the least storage at every $M$ and stays below the dense array throughout, from $0.86$ times its size at $M = 5$ to $0.45$ times at $M = 30$. The global and per-cell layouts need more than the dense array at $M = 5$ ($1.6$ and $1.1$ times), about as much at $M = 10$, and less for $M \ge 20$. The runtime grows with $M$ for all layouts, in the same order as in Fig.~\ref{fig:scaling_nx}.

\begin{figure}[htbp]
\centering
\includegraphics[width=0.625\textwidth]{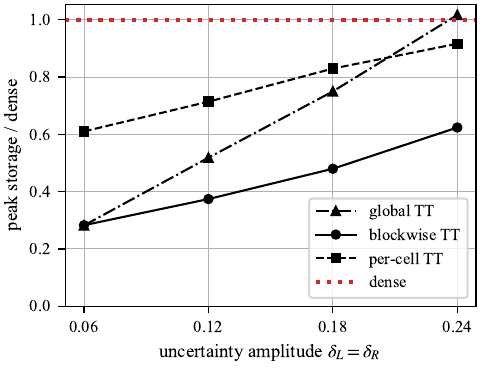} 
\caption{Peak storage relative to the dense array versus the uncertainty amplitude $\delta_L = \delta_R$. $N = 100$, $M = 10$, fixed number of steps.}
\label{fig:storage_uncertainty}
\end{figure}

As the uncertainty amplitude $\delta_L = \delta_R$ increases, the band of uncertain diaphragm positions widens, and the storage of all layouts grows (Fig.~\ref{fig:storage_uncertainty}). For the smallest amplitude, the blockwise and global TT need the same storage, $0.28$ times that of the dense array, and the per-cell TT more than twice as much, $0.61$ times. The global TT grows fastest and slightly exceeds the dense storage at the largest amplitude, where the per-cell TT reaches $0.92$ times, while the blockwise TT grows most slowly and remains well below it, at $0.62$ times.

The scaling studies reflect the early phase of the run. Over the full run to $t = 0.3003$ with $N = 100$ and $M = 10$ (Fig.~\ref{fig:rank_layout}), the blockwise TT has the smallest peak storage, and the shortest total runtime. At every step, it stores at least $26\%$ fewer entries than the per-cell TT, and it never exceeds the dense storage, whereas the per-cell TT exceeds it from $t \approx 0.05$ on and reaches $1.3$ times the dense storage by the last output time. Only from $t \approx 0.23$ on, when the reflected waves have spread the uncertainty over most of the domain and nearly all blocks have the minimum size, does the global TT store fewer entries at the end of a step than the blockwise TT.

\subsection{\texorpdfstring{Six uncertain parameters}{6 uncertain parameters}}
\label{sec:results_param6}
With two parameters, the rank of the link between the two parameter modes reaches its largest possible value $M + 4$, so the savings relative to the dense array are bounded. With more parameters, the links between parameter cores are no longer bounded by $P$ and can be compressed. We therefore apply blockwise TT to a problem with six uncertain parameters, comparing it to an MC reference constructed from $10^{6}$ samples.
\begin{figure}[htbp]
  \centering
  \includegraphics[width=\textwidth]{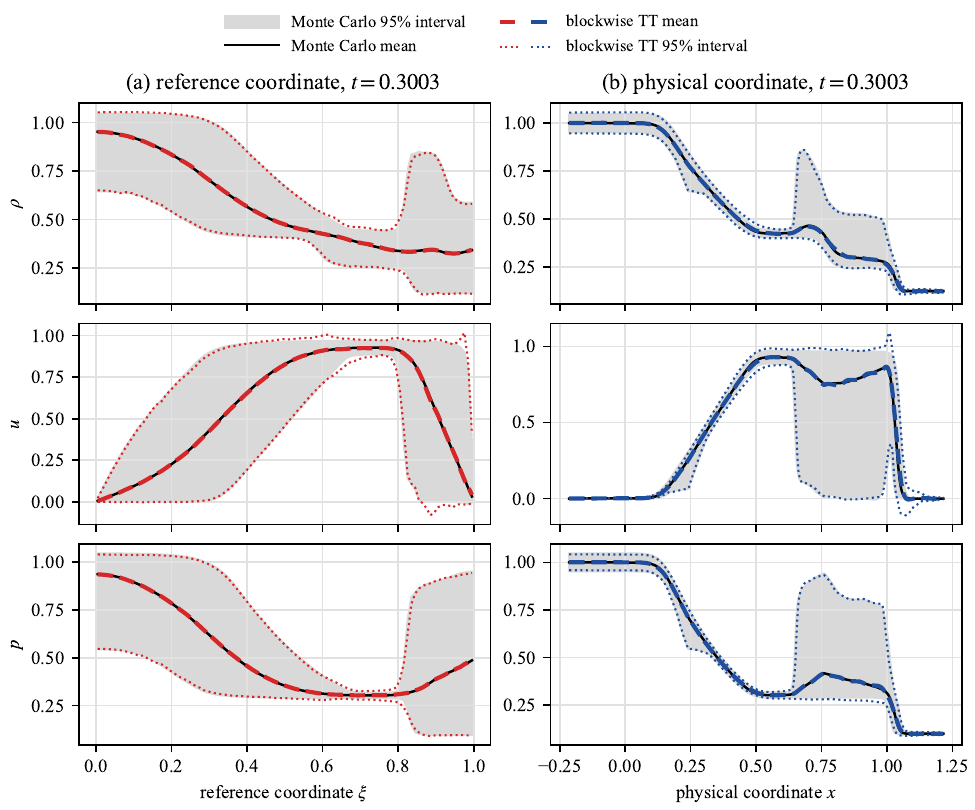}
  \caption{Six uncertain parameters: walls $x_L \in [-\delta_L, \delta_L]$ and
    $x_R \in [1-\delta_R, 1+\delta_R]$ with $\delta_L = \delta_R = 0.24$, and the left and
    right states $\rho_L, p_L, \rho_R, p_R$, each uniform within $\pm 5\%$ of its Sod value;
    the diaphragm is fixed at $x = 1/2$. Mean and central $95\%$ interval, between the $2.5\%$ and $97.5\%$ percentiles for $g=\rho$, $u$, and $p$ (rows) at $t=0.3003$ over (a) the reference coordinate $\xi$ and (b) the physical coordinate $x$. Blockwise TT mean (red/blue dashed) and intervals (red/blue dotted) are obtained via linear interpolation of the Gauss-point reconstruction. The MC reference solution mean is represented by a solid black line and its interval is represented by gray bands. $N=100$, $M=10$, and $\varepsilon_{tt}=10^{-3}$.}
  \label{fig:geo6-statistics}
\end{figure}
To test the method in a higher-dimensional parameter space, we add uncertainty in the initial
states to the geometric uncertainty of Section~\ref{sec:geo_unc}. The walls are uniform,
$x_L \sim \mathcal{U}[-0.24, 0.24]$ and $x_R \sim \mathcal{U}[0.76, 1.24]$
($\delta_L = \delta_R = 0.24$ as in Section~\ref{sec:results_setup}), and the density and pressure on each side of
the diaphragm are independent and uniform within $\pm 5\%$ of their Sod values,
$\rho_L, p_L \sim \mathcal{U}[0.95, 1.05]$, $\rho_R \sim \mathcal{U}[0.11875, 0.13125]$, and
$p_R \sim \mathcal{U}[0.095, 0.105]$. The gas is initially at rest, $\gamma = 1.4$, and the
diaphragm stays fixed at $x = 1/2$. This gives $d = 6$ parameters,
$\mathbf{y} = (x_L, x_R, \rho_L, p_L, \rho_R, p_R)$. Each TT therefore has seven cores: the
physical core followed by one core per parameter in this order, with the two geometric
parameters next to the physical mode. As in Section~\ref{sec:ic_bc}, the initial condition is assembled
exactly in TT format. The wall coordinates couple to the spatial mode through the mapped
diaphragm position $\xi_d(\mathbf{y})$, while each state is constant on its side of the
diaphragm and enters with rank one.

The discretization is that of Section~5.1: $N = 100$ reference cells, $M = 10$ parameter cells
per direction with two padding cells on each side ($P = M + 4 = 14$), WENO3 in $\xi$ and in the
parameters, the Lax--Friedrichs flux, the positivity limiter, and SSP-RK3 with the fixed step
$\Delta t = 0.3\,\Delta\xi\, s_{\min}/a_{\max} = 4.8\cdot10^{-4}$ ($a_{\max} = 2.5$,
$s_{\min} = 0.4$). We take 627 steps to $t = 0.3003$ with the mean-conservative scheme of
Sections~\ref{sec:tt_cons_round} and~\ref{sec:block_flux_share}. The TT tolerance is $\varepsilon_{\rm tt} = 10^{-3}$, with rank cap
$r_{\max} = 128$, which is never reached. Here we use $\varepsilon_{\rm tt} = 10^{-3}$ instead of the $\varepsilon_{\rm tt} = 10^{-4}$ tolerance of Section~\ref{sec:results_setup}, since the cost grows quickly with the number of parameters. Each cross approximation stops once its relative
$L^2$ error on a fixed set of 1000 random entries of the approximated tensor falls below
$\varepsilon_{\rm tt}$, or once that error stops decreasing; the truncation tolerance inside
the two-site sweep is $\varepsilon_{\rm tt}/10$.

\begin{figure}[htbp]
  \centering
  \includegraphics[width=\textwidth]{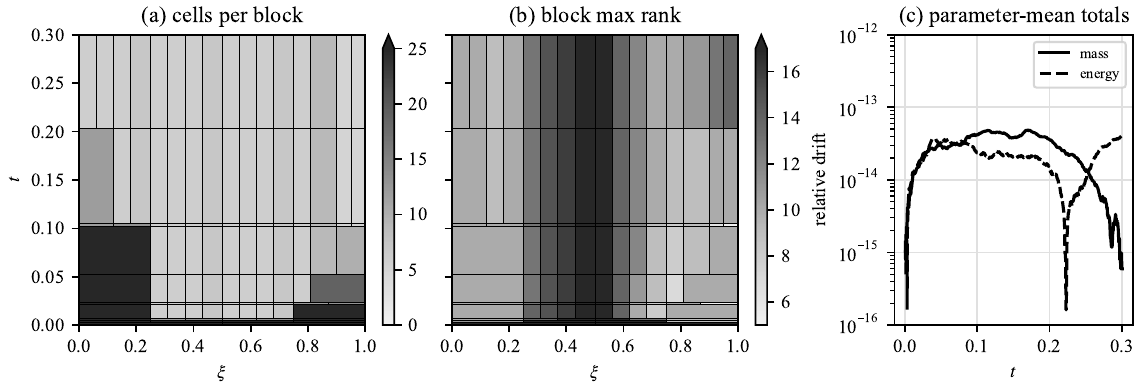}
  \caption{Ranks and block layout for the six-parameter test of
    Fig.~\ref{fig:geo6-statistics}. (a) Cells per block and (b) largest TT rank per block, over reference space and time. (c) Relative drift $|\mathcal{M}_\nu(t) - \mathcal{M}_\nu(0)|/|\mathcal{M}_\nu(0)|$ in time of the mean total mass and energy using the conservative scheme. $N=100$, $M=10$, and $\varepsilon_{tt}=10^{-3}$.}
  \label{fig:geo6-ranks-layout}
\end{figure}

Unlike Section~\ref{sec:results}, the run starts from a single global TT per conserved variable over all
$N = 100$ cells and lets the adaptation of Section~\ref{sec:split_merge} build the block structure. Adaptation
runs after the first step and then every 5 steps with the shock and uncertainty indicators of
Section~\ref{sec:block_indicators} and the thresholds of Eqn.~\ref{eqn:ind_thresholds}. With six parameters the links
between the parameter cores are no longer bounded by $P$, so we take $R = r_{\max}$ for the rank
thresholds. Blocks have between 5 and 20
cells, a merge must keep its error below $\varepsilon_{\rm merge} = 2\cdot10^{-3}$, and up to four
merge sweeps are made per adaptation. The single initial block, of rank 33, is halved twice in
the first five steps. The shock and uncertainty indicators then split the blocks over the band
of uncertain waves down to 6 or 7 cells by $t \approx 0.007$, and later splits follow the waves to
the walls, giving a final layout of 16 blocks of 5 to 9 cells (Fig.~\ref{fig:geo6-ranks-layout}(a)).
 After the start-up the largest ranks are 15, 17, and 15, occurring in the band of uncertain diaphragm positions.(Fig.~\ref{fig:geo6-ranks-layout}(b)).

Fig.~\ref{fig:geo6-statistics} compares the statistics of the reference and blockwise representation at $t = 0.3003$ in both coordinates. In physical coordinates the intervals retain similar structure to the
two-parameter case (Fig.~\ref{fig:profiles_mean_x}). Intervals near the walls are wide since
reflected waves depend on the wall positions, and narrow bands elsewhere arise from the additional $\pm 5\%$
uncertainty in the states. The parameter-mean mass and energy drift by at most
$5\cdot10^{-14}$ (Fig.~\ref{fig:geo6-ranks-layout}(c)). The run takes $5.1$\,hrs on 8 threads, and the
blockwise TT stores at most $8.0\cdot10^4$ entries, against $3N P^6 \approx 2.3\cdot10^9$ for the
dense array. In the parameter direction, the ranks of the largest block decrease from the
geometric to the state modes, so the four state parameters add little storage.

\begin{figure}[htbp]
  \centering
\includegraphics[width=\textwidth]{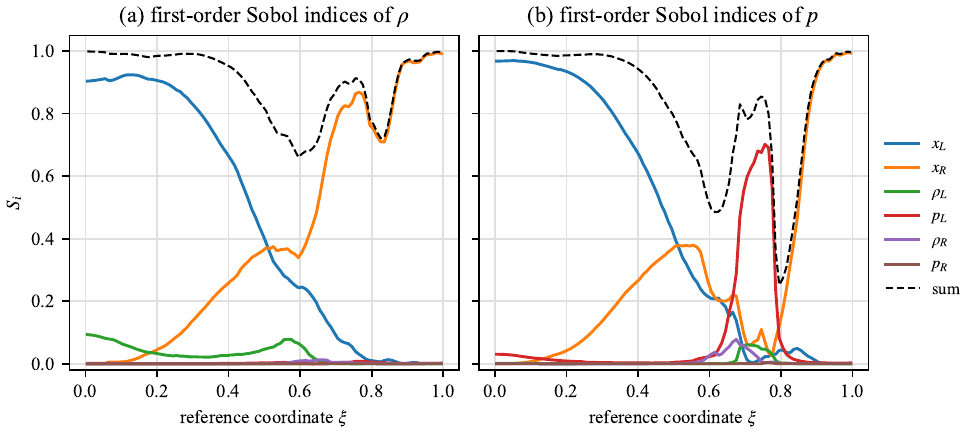}
  \caption{First-order Sobol indices $S_i = \mathrm{Var}(\mathbb{E}[g \mid y_i]) / \mathrm{Var}(g)$ of (a) $g=\rho$ and (b) $g=p$ over reference coordinate space at $t=0.3003$ for the six-parameter case of Fig.~\ref{fig:geo6-statistics}. The dashed line represents the sum of the first order indices; where it is below one, the parameters interact. $N=100$, $M=10$, and $\varepsilon_{tt}=10^{-3}$.}
  \label{fig:geo6-sensitivity}
\end{figure}

The blockwise TT representation can also be used to extract variance-based sensitivities at a negligible cost, as the first-order Sobol
index $S_i = \mathrm{Var}(\mathbb{E}[g \mid y_i]) / \mathrm{Var}(g)$ follows from contracting the
train over all parameters but $y_i$. Fig.~\ref{fig:geo6-sensitivity} shows these indices for
$\rho$ and $p$. Near each wall the variance is almost entirely due to that wall. In the pressure
plateau between the contact and the reflected shock, $p$ is controlled mainly by $p_L$, while the
wave positions there barely matter. Where the waves of the realizations overlap, the sum of the
first-order indices drops to about $0.66$ for $\rho$ and $0.25$ for $p$, indicating that the parameters
interact.

%
%
\section{\texorpdfstring{Discussion}{Discussion}}
\label{sec:discussion}
The global TT is limited by its single leading rank, which the few cells containing the uncertain waves set for the whole domain, and the per-cell TT by the replication of the parameter cores in every cell and its large number of TT operations. The blockwise TT shares the parameter cores within a block and bounds the leading rank by the block size, so a discontinuity raises the ranks only of the blocks it passes through. Since cross approximation dominates the runtime, the cost of a layout is essentially set by the number of cross approximations and the ranks of the trains they act on, and the blockwise TT keeps both small. With two parameters, the parameter link is at full rank in the blockwise and per-cell layouts, so the compression comes from the link between the physical and the parameter modes, which is exactly the link the partition controls. The savings relative to the dense array are therefore bounded by a factor of a few, while those relative to the other two layouts come entirely from the partition.

A fixed partition can only bound the leading rank by its block size everywhere, whereas the adaptive partition chooses the block size locally, small where the uncertain waves are and large elsewhere (Fig.~\ref{fig:rank_layout}). Once the reflected waves have spread the uncertainty over most of the domain, the partition approaches the minimum block size everywhere, and the global TT stores fewer entries than the blockwise TT, although over the whole run the blockwise TT keeps the smaller peak storage and runtime (Section~\ref{sec:results_cost}).

For the range of $M$ considered, the accuracy of the statistics is limited by the stochastic finite volume discretization rather than by the TT format (Fig.~\ref{fig:convergence}). The slightly larger blockwise and global errors of the mean and the standard deviation at $M = 40$ are consistent with a truncation error that is relative to the norm of the whole train. Conservation of the parameter-mean mass and energy holds for every layout and partition, since the correction of Eqn.~\ref{eqn:cons_correction} adds no rank and the interface fluxes are single-valued (Fig.~\ref{fig:conservation}(a, b)), and, unlike with standard rounding, it is independent of the TT tolerance (Fig.~\ref{fig:conservation}(c)). It therefore need not be bought with a tighter tolerance and the higher ranks that would entail. With six uncertain parameters (Section~\ref{sec:results_param6}), the parameter links are no longer at full rank, and the blockwise TT stores more than four orders of magnitude fewer entries than the
dense array. The ranks decrease from the geometric to the state modes, so the four state parameters add little storage. The trains can also be leveraged to give first-order Sobol indices by contraction, without further solves.

%
%
\section{\texorpdfstring{Conclusions}{Conclusions}}
\label{sec:conclusions}
We have presented an adaptive blockwise TT representation for stochastic finite volume methods with uncertain geometry. The SFV solution is partitioned along its physical mode into blocks, each stored as a separate TT for every conserved variable, and the partition is adapted by split and merge operations driven by shock, uncertainty, and rank indicators. A rounding that restores the parameter-mean totals without adding rank, together with single-valued interface fluxes, conserves the parameter-mean mass and energy up to round-off for any partition and independently of the TT tolerance.

On a Sod shock tube with uncertain wall positions, the blockwise TT reproduced the mean, and the percentile intervals of a dense reference in reference and in physical coordinates as accurately as the global and per-cell TT. It was the most compact layout over all mesh sizes, parameter resolutions, and uncertainty amplitudes in the scaling studies, with a peak storage $41\%$ below that of the global TT and $38\%$ below that of the per-cell TT on the finest mesh. It was also the fastest, $2.2$ times faster than the global TT and $5.1$ times faster than the per-cell TT on the finest mesh, with about seven times fewer trains than the per-cell TT.

These advantages come from adapting the partition to where the uncertainty is. A single global TT lets the few cells containing the uncertain waves set the ranks of the whole domain, and a separate TT per cell loses the parametric structure that neighboring cells share. The blockwise TT keeps the parameter cores shared within each block and confines rank growth to the blocks the waves pass through, so it combines the strengths of both layouts at a lower cost than either. Its conservation of the parameter-mean mass and energy does not depend on the TT tolerance, so it need not be paid for with a tighter tolerance and higher ranks. The gains are largest when the uncertainty is localized: for the smallest uncertainty amplitude, the blockwise and global TT needed about the same storage, and once the uncertain waves had filled the domain, the global TT became more compact. Since cross approximation dominates the cost of every layout, the runtime follows the number of cross approximations and the ranks they act on, and the per-cell TT, which performs the most operations, is the slowest layout. With two parameters the savings relative to the dense array are bounded by a factor of a few, because the link between the two parameter modes is at full rank. With six uncertain parameters, the blockwise TT closely matched a Monte Carlo reference with
$10^{6}$ samples, and stored about $3\times10^{4}$ times fewer entries than the dense array. Problems with even more uncertain parameters, more physical dimensions, and longer simulations, in which the parameter cores themselves are compressible and the dense array grows fastest, are where the blockwise representation is expected to pay off most, and are the subject of future work.

\section*{\label{sec:ack} Acknowledgments}
Research presented in this article was supported by the Laboratory Directed Research and Development program of Los Alamos National Laboratory under project number LDRD-20250032DR.  The Los Alamos unlimited release number is LA-UR-26-30606.

During the preparation of this work, the authors used Claude Opus 5.5 (Anthropic), to assist with portions of the code development. The tool was also used to assist with drafting portions of the manuscript text. All AI-assisted code was reviewed, tested, and validated by the authors, and all AI-assisted text was reviewed and revised by the authors. The authors take full responsibility for the content of this publication.

\bibliographystyle{siamplain}
\bibliography{references}

\end{document}

%% file: ex_shared.tex
\usepackage{lipsum}
\usepackage{amsfonts}
\usepackage{graphicx}
\usepackage{algorithm}
\usepackage{bbm}
\usepackage{algpseudocode}
\usepackage{epstopdf}
\usepackage{microtype}
\usepackage{placeins}
\ifpdf
  \DeclareGraphicsExtensions{.eps,.pdf,.png,.jpg}
\else
  \DeclareGraphicsExtensions{.eps}
\fi

\newsiamremark{assumption}{Assumption}

\title{Adaptive Blockwise Tensor Trains for Stochastic Finite-Volume Methods with Uncertain Geometry\thanks{Submitted to the editors DATE.
\funding{The Laboratory Directed Research and Development (LDRD) program of Los Alamos National Laboratory supported the work of WG, SW, SA, ST, and WT under project number 20250032DR. Los Alamos National Laboratory is operated by Triad National Security, LLC, for the National Nuclear Security Administration of U.S. Department of Energy (Contract No. 89233218CNA000001).}}}

\author{
  William P. Gammel \thanks{XCP-8, Los Alamos National Laboratory, Los Alamos, NM, USA 
  (\email{wpgam@lanl.gov}).}
\and Steven T. Walton \thanks{T-5, Los Alamos National Laboratory, Los Alamos, NM, USA 
  (\email{stevenw@lanl.gov}).} 
\and Stephen A. Andrews \footnotemark[3]
\and Svetlana Tokareva \thanks{T-5, Los Alamos National Laboratory, Los Alamos, NM, USA 
  (\email{tokareva@lanl.gov}).}  
\and William T. Taitano \thanks{T-5, Los Alamos National Laboratory, Los Alamos, NM, USA 
  (\email{taitano@lanl.gov}).}  
}

\usepackage{amsopn}
